\documentclass[12pt]{amsart}
\usepackage{amssymb,amsthm}
\usepackage[all]{xy}
\theoremstyle{plain}
\newtheorem{thm}{Theorem}[section]
\newtheorem{lem}[thm]{Lemma}
\newtheorem{prop}[thm]{Proposition}
\newtheorem{cor}[thm]{Corollary}

\theoremstyle{remark}
\newtheorem{rem}[thm]{Remark}
\newtheorem{rems}[thm]{Remarks}

\theoremstyle{definition}
\newtheorem{defn}[thm]{Definition}

\numberwithin{equation}{thm}

\DeclareMathOperator{\Ann}{Ann}
\DeclareMathOperator{\HH}{HH}
\DeclareMathOperator{\coh}{H}
\DeclareMathOperator{\Hom}{Hom}
\DeclareMathOperator{\PHom}{PHom}
\DeclareMathOperator{\Ext}{Ext}

\DeclareMathOperator{\Ker}{Ker}
\DeclareMathOperator{\Ima}{Im}
\DeclareMathOperator{\rad}{rad}

\DeclareMathOperator{\Span}{Span}

\DeclareMathOperator{\id}{id} 

\DeclareMathOperator{\soc}{soc}

\newcommand{\ul}{\underline{\lambda}}

\newcommand{\Z}{\mathbb{Z}}
\newcommand{\Zl}{\Z/\ell\Z}

\newcommand{\ot}{\otimes}
\newcommand{\Proj}{\text{Proj}\,}

\newcommand{\bu}{\bullet}
\newcommand{\ula}{\underline\lambda}
\newcommand{\va}{V_A^c}
\DeclareMathOperator{\Spec}{Spec}

\DeclareMathOperator{\kt}{k\langle \tau_{\ul} (t)\rangle}

\title[Varieties for Modules]
{Varieties for Modules of Quantum Elementary Abelian Groups}
\author{Julia Pevtsova}
\address{Department of Mathematics, University of Washington, Seattle, WA
98195, USA}
\email{julia@math.washington.edu}
\thanks{The first author was supported by
NSF grant \#{}DMS--0500946}
\author{Sarah Witherspoon}
\address{Department of Mathematics, Texas A\&M University, College Station,
TX 77843, USA}
\email{sjw@math.tamu.edu}

\thanks{The second author was supported by the Alexander von
Humboldt Foundation and by
NSF grants \#{}DMS--0422506 and \#{}DMS-0443476.}
\date{March 14, 2006}

\begin{document}

\begin{abstract}
We define a rank variety for a module of a noncocommutative Hopf algebra 
$A = \Lambda \rtimes G$ where 
$\Lambda = k[X_1, \dots, X_m]/(X_1^{\ell}, \dots, X_m^{\ell})$, $G = (\Zl)^m$, 
and $\text{char } k$ does not divide $\ell$, 
in terms of certain subalgebras of 
$A$ playing the role of ``cyclic shifted subgroups".  
We show that the rank variety of a finitely generated module $M$ 
is homeomorphic to the support variety of $M$ defined in terms of the 
action of the cohomology algebra of $A$.
As an application we derive a theory of rank varieties for the algebra 
$\Lambda$. When $\ell=2$,
rank varieties for $\Lambda$-modules were constructed by Erdmann and Holloway 
using the representation theory of the Clifford algebra. 
We show that the rank varieties we obtain for $\Lambda$-modules coincide with those of
Erdmann and Holloway.  
\end{abstract}

\maketitle

\section{Introduction}

The theory of varieties for modules of a finite group $G$ began with
the groundbreaking work of Quillen \cite{Q}, a stratification of the maximal
ideal spectrum of the cohomology ring of $G$ into pieces indexed by
elementary abelian subgroups. This idea was taken  further by
Avrunin and Scott \cite{AS}, to a stratification of an affine variety
associated to any finitely generated module. These results depended
on earlier work of Venkov \cite{Venk} and Evens \cite{Ev}, 
showing that the cohomology of $G$,
a graded commutative ring, is finitely generated.

The theory took a different twist with the introduction by Carlson \cite{C1}
of the {\it rank variety} for a module of an elementary abelian group $E$. 
The rank variety is yet another geometric  invariant of a module, and is 
defined in terms of {\it cyclic shifted subgroups} of $E$. 
Carlson conjectured that the variety arising from the action of cohomology, 
and the rank variety defined purely in terms of representation-theoretic 
properties of a module, coincide. The conjecture was 
proven by Avrunin and Scott \cite{AS}. 

This theory was adapted to restricted Lie algebras by Friedlander and 
Parshall \cite{FP}. It was then further generalized to other finite group 
schemes (see \cite{FPv1,SFB1,SFB2}) 
based upon the fundamental theorem of Friedlander and Suslin stating
that the cohomology of any finite group scheme, 
or equivalently finite dimensional cocommutative Hopf algebra, is
finitely generated \cite{FS}.
In particular, the notion of rank variety was recently 
generalized to all finite group schemes by Friedlander and the first author
\cite{FPv1}. 
One important aspect of the rank variety in the context of finite group 
schemes is that it satisfies the ultimate generalization of the Avrunin-Scott 
Theorem:  The rank variety of a module defined in a purely 
representation-theoretic way is homeomorphic 
to the support variety defined cohomologically. The interplay between the two seemingly very different 
descriptions of the variety of a module  allows for applications both in cohomology and in representation theory.

Much less is known in the context of
finite dimensional noncocommutative Hopf algebras.
Ginzburg and Kumar computed the cohomology rings of quantum groups at
roots of unity, and these happen to be finitely generated \cite{GK}.
This fact allowed mathematicians to start development of the 
theory of support varieties for modules of these small quantum groups 
(see \cite{Ost}, \cite{PW}).
However it appears difficult to give an equivalent representation-theoretic
definition of variety for these quantum groups in general. Even less
has been done for other types of finite dimensional noncocommutative
Hopf algebras, and in particular it is an open question as to whether
their cohomology is finitely generated.

In this paper, we have very modest goals. We only consider Hopf
algebras that are quantum analogues of elementary abelian groups, 
namely tensor products of Taft algebras (which are also Borel subalgebras
of $u_q(sl_2^{\times m})$). We define the rank variety of a module for such a Hopf
algebra (Definition \ref{defn:rank}), 
giving the first definition of rank varieties for modules of
a noncocommutative
Hopf algebra. The cohomology of a tensor product of Taft algebras is
finitely generated, so we may also associate a
support variety, defined cohomologically, to any module (\ref{eqn:suppMN}). 
We show that the rank variety of any
finitely generated module is homeomorphic to the support variety
(Theorem \ref{thm:equiv}), thus 
providing an analogue of the Avrunin-Scott Theorem in our context. 
We use ``Carlson's modules" $L_\zeta$ as our main tool and apply the techniques developed  
in \cite{EH2} and \cite{EHSST} in the study of support varieties defined via Hochschild cohomology.
We expect that our results will shed light on
the problem of constructing a rank variety for a broader class of 
finite dimensional Hopf algebras, including the small quantum groups.

One of the most important applications of the identification of the 
rank and support varieties in the setting of finite group schemes is the proof 
of the ``tensor product property"   which expresses the variety of a tensor 
product as the intersection of varieties (see \cite{AS}, \cite{FP}, \cite{FPv1}, \cite{SFB2}).
Another common application is a classification 
of thick tensor ideal subcategories in the stable module  category (see \cite{BCR2}, \cite{FPv2}).
Both of these applications will be addressed in a sequel to this paper.

Our results have consequences beyond Hopf algebras. A tensor product
of Taft algebras is isomorphic to a
skew group algebra $A=\Lambda\rtimes G$ where the
group $G \cong (\Zl)^{m}$ is elementary abelian (in nondefining characteristic)
and $\Lambda = k[X_1, \dots, X_m]/(X_1^{\ell}, \dots, X_m^{\ell})$. 
When $\ell =2$, that is the generators of $\Lambda$ have square $0$,
Erdmann and Holloway have used Hochschild cohomology to define support
varieties for $\Lambda$-modules \cite{EH2}, applying a theory of varieties for
modules of algebras initiated by Snashall and Solberg \cite{SS}.
The support variety of a $\Lambda$-module in this
case is equivalent to a rank variety defined representation-theoretically
by Erdmann and Holloway.
Their approach is quite different from ours: 
They use a ``stable map description" of the rank variety  and
representation theory of the Clifford algebra. 
In this paper we use the extension of $\Lambda$ to $A$ to give definitions
of support and rank varieties for $\Lambda$-modules more generally (see
(\ref{eqn:ranklambda}) and (\ref{eqn:supplambda})), 
that is for {\em any} $\ell$ not divisible by the characteristic of the 
field $k$, and to show that the varieties we obtain are homeomorphic
(Corollary \ref{cor:equivlambda}). 
In case the generators of $\Lambda$ have square 0, our
varieties coincide with those of Erdmann and Holloway, giving
an alternative approach to their theory. 
In order to make this connection, we found it necessary to record
some basic facts relating cohomology and Hochschild cohomology of
finite dimensional Hopf algebras in an appendix.

When this article was nearly
complete, the authors learned that Benson, Erdmann, and Holloway
had found a different way to define rank varieties for $\Lambda$-modules
for arbitrary $\ell$, involving an algebra extension of $\Lambda$ that is a tensor
product of $\Lambda$ with a twisted group algebra of $G$ \cite{BEH}. 
Their algebra
extension has some features in common with ours, leading to a parallel
theory. We thank Benson, Erdmann, and Holloway  for some very helpful
conversations.

We thank \'Ecole Polytechnique F\'ed\'erale de Lausanne and
Universit\"at M\"unchen for their 
hospitality during the preparation of this paper. 

\vspace{.1in}

Throughout this paper, $k$ will denote
an algebraically closed field containing a primitive $\ell$th root of unity $q$.
All tensor products and dimensions will be over $k$ unless otherwise indicated.
We shall use
the notation $V^\#$ for the $k$-linear dual of a finite dimensional
vector space $V$.

\section{Quantum analogues of cyclic shifted subgroups}

Let $m$ be a positive integer and let
$G$ denote the group $(\Z/\ell\Z)^m$ with generators $g_1,\ldots,g_m$.
Define an action of $G$ by automorphisms on the polynomial ring $R=k[X_1,\ldots,X_m]$ by
setting
$$
  g_i\cdot X_j = q^{\delta_{ij}}X_j
$$
for all $i,j$, where $\delta_{ij}$ is the Kronecker delta. 
Let $\widetilde{A}=R\rtimes G$, the skew group algebra,
that is $\widetilde{A}$ is a free left
$R$-module having $R$-basis $G$, with the 
semidirect (or smash) product multiplication
$$
  (rg) (sh) = r (g\cdot s) gh
$$
for all $r,s\in R$ and $g,h\in G$.
Then $\widetilde{A}$ is a Hopf algebra with 
$$
  \Delta(X_i)=X_i\ot 1 + g_i\ot X_i, \hspace{.3in}\Delta(g_i)=g_i\ot g_i,
$$
$\varepsilon(X_i)=0$, $\varepsilon(g_i)=1$,
$S(X_i)=-g_i^{-1}X_i$, and $S(g_i) = g_i^{-1}$, for all $i$. 
Letting $h_1=1$ and $h_j=\prod _{i=1}
^{j-1} g_i \ (2\leq j\leq m)$, we have
\begin{equation}\label{eqn:q-commute}
  X_jh_j\cdot X_ih_i = q X_ih_i\cdot X_jh_j \ \ \mbox{ for all }j>i.
\end{equation}
The following consequence of this $q$-commutativity of the elements $X_ih_i$ 
will be most essential in what follows.

\begin{lem}\label{lem:tau}
For any $\lambda_1,\ldots,\lambda_m \in k$,
$\displaystyle{\left(\sum_{j=1}^m\lambda_j X_jh_j\right)^{\ell}=
  \sum_{j=1}^m \lambda_j^{\ell} X_j^{\ell}}$.
\end{lem}

\begin{proof}
This is a consequence of the $q$-binomial formula which in this context gives,
for all $n\leq \ell$ and $j> i$,
$$
  (\lambda_i X_i h_i + \lambda_jX_j h_j)^n =\sum_{s=0}^n
   \frac{(n)_q!}{(s)_q! (n-s)_q!}
  (\lambda_iX_i h_i)^s(\lambda_jX_j h_j)^{n-s},
$$
where 
$(s)_q = 1+q+q^2+\cdots +q^{s -1}$, 
$(s)_q ! = (s)_q(s-1)_q\cdots (1)_q$, 
and $(0)_q!=1$ by definition.
If $n=\ell$, the coefficients of $\lambda_i^{\ell}X_i^{\ell}$ and
$\lambda_j^{\ell}X_j^{\ell}$ should be interpreted to be $1$.
As $q$ is a primitive $\ell$th root of 1, induction on $m$ yields the
desired result.
\end{proof}

Another application of the $q$-binomial formula, to $\Delta(X_i^{\ell})=
(X_i\otimes 1 + g_i\otimes X_i)^{\ell}$, shows that 
the ideal $(X_1^{\ell},\ldots,X_m^{\ell})$ is a Hopf ideal.  
Thus
$$
  A = \widetilde{A}/(X_1^{\ell},\ldots,X_m^{\ell})
$$
is a Hopf algebra of dimension $\ell^{2m}$, a tensor product of $m$ 
copies of a Taft algebra, a quantum analogue of an elementary abelian group.
We may identify $A$ with
the skew group algebra $\Lambda\rtimes G$ where 
$$
  \Lambda = k[X_1,\ldots,X_m]/(X_1^{\ell},\ldots,X_m^{\ell}),
$$
a truncated polynomial algebra. We will primarily be interested in
the finite dimensional Hopf algebra $A$ in this paper, but will need to use
$\widetilde{A}$ as well in some of the proofs.
Note that since $A$ is a finite dimensional Hopf algebra, it is a 
Frobenius algebra \cite[Thm.\ 2.1.3]{montgomery93}, and in particular
is self-injective.

We now introduce algebra maps $\bf \tau_{\ul}$ 
which will play the role of ``cyclic shifted subgroups" (see \cite[II]{Ben}) 
or $p$-points (\cite{FPv1}) for the algebra $A$. 
By Lemma \ref{lem:tau}, for each point $\ul = [\lambda_1:\ldots: \lambda_m]$ in
$k$-projective space
$\mathbb P^{m-1}$,
there is an embedding of algebras
\begin{equation}\label{defn:taulambda}
  \tau_{\ul} : k[t]/(t^{\ell}) \rightarrow A
\end{equation}
defined by $\tau_{\ul}(t)=\sum_{i=1}^m \lambda_i X_i h_i$.
Denote the image of $\tau_{\ul}$ by $k\langle \tau_{\ul}(t)\rangle$.

\begin{lem}\label{lem:free}
Let   $\ul = [\lambda_1:\ldots: \lambda_m] \in \mathbb P^{m-1}$.
Then $A$ is free as a left (respectively, right)
$k\langle \tau_{\ul} (t)\rangle$-module, with $k\langle\tau_{\ul}(t)\rangle$-basis
$$
 {\mathcal B}=  \{X_2^{a_2}\cdots X_{m}^{a_{m}} g_1^{b_1}
  \cdots g_m^{b_m}\mid 0\leq a_i,b_i\leq \ell -1\}
$$
in case $\lambda_1\neq 0$. Analogous statements hold if $\lambda_i\neq 0$
for other values of $i$.
\end{lem}

\begin{proof}
The statement that $A$ is free over $k\langle \tau_{\ul}(t)\rangle$ follows
from a general result of Masuoka on coideal subalgebras \cite{Mas}.
We shall use the explicit basis $\mathcal B$ however, and so we provide
a proof for completeness.
 We shall prove that $\mathcal B$ is a free $k\langle\tau_{\ul}(t)\rangle$-basis
of $A$ as a left $k\langle\tau_{\ul}(t )\rangle$-module. 
That it is also a basis of $A$ as a right module is proved similarly.

We may assume that $\lambda_1 = 1$.
 Since the number of elements in $\mathcal B$ is $\ell^{2m-1}=
\dim_k A/ \dim_k k\langle\tau_{\ul}(t )\rangle$, it suffices to show that $A =
k\langle\tau_{\ul}(t )\rangle \mathcal B$.
\sloppy
{

}

We use induction on $a_1$ to show that  $X_1^{a_1}\cdots X_m^{a_m} g_1^{b_1}
  \cdots g_m^{b_m} \in k\langle\tau_{\ul}(t )\rangle \mathcal B$ for any
choice of exponents $0\leq a_i,b_i\leq \ell -1$. 
The statement is trivial for $a_1 = 0$.   Assume it is proved for all
monomials with $a_1 < n\leq \ell -1$. It remains to show 
that $X_1^{n}X_2^{a_2}\cdots X_{m}^{a_{m}} g_1^{b_1}
  \cdots g_m^{b_m} \in k\langle\tau_{\ul}(t )\rangle \mathcal B$.
The defining relations on $X_i$ and $g_j$, together with the
definition of $\tau_{\ul}(t )$ given above (\ref{defn:taulambda}), immediately imply that 
$$
  X_1^{n}X_2^{a_2}\cdots X_{m}^{a_{m}} g_1^{b_1}
  \cdots g_m^{b_m} - \tau_{\ul}(t ) X_1^{n-1}X_2^{a_2}\cdots
  X_{m}^{a_{m}} g_1^{b_1}\cdots g_m^{b_m}
$$ 
is a sum of monomials $X_1^{a'_1}\cdots
X_m^{a'_m} g_1^{b'_1}
  \cdots g_m^{b'_m}$ for some exponents $a_i',b_i'$ with
$a_1'<n$. The statement follows by induction. 
\end{proof}

To every point $\ul$ in $k$-projective space $\mathbb P^{m-1}$ we associate two special {\it left}
$A$-modules: $V(\ul )$ and 
$V^\prime(\ul )$, which will be used extensively throughout the paper.
 We point out that our modules are different from those used in \cite{EH1}, \cite{EH2}
even though we choose to use similar names for them. 
As will be shown in Corollary~\ref{rankvlambda}, they share one of the 
main properties with the
modules introduced in \cite{EH1}: The rank variety of each of $V(\ul )$ and  
$V^\prime(\ul )$ will be the point $\ul \in \mathbb P^{m-1}$.

\quad

For each $\ul\in {\mathbb{P}}^{m-1}$, let 
$$
  V(\ul)= A\cdot\tau_{\ul}(t)^{\ell -1} \ \ \ \mbox{ and }
  \ \ \ V'(\ul)=A\cdot \tau_{\ul}(t),
$$
that is $V(\ul)$ (respectively, $V'(\ul)$) is the left ideal generated
by $\tau_{\ul}(t)^{\ell -1}$ (respectively, $\tau_{\ul}(t)$).

\vspace{0.1in}
 Recall that for an $A$-module $M$, the {\it Heller shift} of $M$, denoted 
$\Omega(M)$, is the kernel of the projection $P(M) \to M$ where $P(M)$ is the projective cover of $M$. Similarly, $\Omega^{-1}(M)$ is the cokernel of the embedding of $M$ into its injective hull.

\begin{lem}
\label{lem:vlambda}

For each $\ul\in {\mathbb{P}}^{m-1}$ we have:

\begin{itemize}

\item[(i)] $V(\ul)\cong k\uparrow_{k\langle\tau_{\ul}(t)\rangle}^A =
  A\ot_{k\langle\tau_{\ul}(t)\rangle}k$.
  
\item[(ii)] The restriction $V(\ul)\!\downarrow_{k\langle\tau_{\ul}(t)\rangle}$ 
contains the trivial module as a direct summand. 
In particular, $V(\ul)$ is not projective as a $k\langle\tau_{\ul}(t)\rangle$-module.

\item[(iii)] $\dim_k V(\ul) =\ell ^{2m-1}$, $\dim_k V^\prime(\ul) =(\ell
-1)\ell^{2m-1}$, and there is a short exact sequence of $A$-modules
$
  0\rightarrow V'(\ul)\stackrel{\iota}{\longrightarrow} 
  A\stackrel{\pi}{\longrightarrow} V(\ul)\rightarrow 0.
$

\item[(iv)]
$\hspace{.2in}\displaystyle{
  \cdots\stackrel{\cdot\tau_{\ul}(t)}{
  \relbar\joinrel\longrightarrow} A 
  \stackrel{\cdot\tau_{\ul}(t)^{\ell -1}}{\relbar\joinrel\relbar\joinrel
  \relbar\joinrel\longrightarrow} 
  A \stackrel{\cdot\tau_{\ul}(t)}
  {\relbar\joinrel\longrightarrow}
  A \stackrel{\cdot\tau_{\ul}(t)^{\ell -1}}{\relbar\joinrel\relbar\joinrel
  \relbar\joinrel\longrightarrow} 
  V(\ul)\rightarrow 0}
$

is a minimal projective resolution of $V(\ul)$, and
$$
  0\rightarrow V(\ul)\rightarrow A \stackrel
  {\cdot\tau_{\ul}(t)}{\relbar\joinrel\longrightarrow} A 
  \stackrel{\cdot\tau_{\ul}(t)^{\ell -1}}{\relbar\joinrel\relbar\joinrel
  \relbar\joinrel\longrightarrow} 
  A \stackrel{\cdot\tau_{\ul}(t)}
  {\relbar\joinrel\longrightarrow}
  A \stackrel{\cdot\tau_{\ul}(t)^{\ell -1}}{\relbar\joinrel\relbar\joinrel
  \relbar\joinrel\longrightarrow} \cdots
$$
is a minimal injective resolution of $V(\ul)$. 
This remains true if 
$V(\ul)$ is replaced by $V'(\ul)$, with appropriate changes in the powers
of $\tau_{\ul}(t)$.

\item[(v)]
$\Omega (V(\ul))\cong V'(\ul)$ and $\Omega(V'(\ul))
\cong V(\ul)$.

\end{itemize}

\end{lem}

\begin{proof} 
(i) Define a map $\phi:A\times k\rightarrow V(\ul) = A\tau_{\ul}(t)^{\ell -1}$ 
by $\phi(a,c) = ca\tau_{\ul}(t)^{\ell -1}$ for all $a\in A,c\in k$,
a $k$-bilinear map that commutes with left multiplication by elements of $A$.
Note that $\phi(a  \tau_{\ul}(t) , c)=0=\phi(a,\tau_{\ul}(t)\cdot c)$,
the latter equality due to the trivial action of $k\langle \tau_{\ul}(t)
\rangle $ on $k$.
Thus $\phi$ induces an $A$-map from the tensor product 
$A\ot_{k\langle\tau_{\ul}(t)\rangle} k$ to $V(\ul)$.  One readily checks that this map gives a bijection 
between the $k$-bases ${\mathcal B}\ot 1$ of 
$A\ot_{k\langle\tau_{\ul}(t)\rangle} k$ and ${\mathcal B}\tau_{\ul}(t)^{\ell -1}$
of $V(\ul)$ where ${\mathcal B}$ is defined in Lemma \ref{lem:free}.

(ii) The trivial $k\langle\tau_{\ul}(t)\rangle$-submodule 
$1\ot k$ of $A\ot_{k\langle\tau_{\ul}(t)\rangle} k$
is complemented by the $k$-linear span of $({\mathcal B} - \{1\})\ot
_{k\langle\tau_{\ul}(t)\rangle} k$.

(iii) By the proof of (i) above,
$\mathcal B$ is in bijection with a $k$-basis
of $V(\ul)$, so $\dim_k V(\ul)=\ell^{2m-1}$. 
Similarly, a $k$-basis of $V'(\ul)$ is $\cup_{i=1}^{\ell -1} {\mathcal B}
\tau_{\ul}(t)^i$, of cardinality $(\ell -1)\ell^{2m-1}$.
The map $\iota : V'(\ul)\rightarrow A$ in the statement of the lemma
is inclusion, and the map $\pi : A\rightarrow V(\ul)$ is given by
$\pi(a)=a\ot 1\in A\ot_{k\langle\tau_{\ul}(t)\rangle} k\cong V(\ul)$.
Again by considering bases of each of these modules, the sequence 
given in the lemma is seen to be exact.

(iv)
Note that the Jacobson radical of $A$ is $\rad(A)= A\cdot \rad (\Lambda)$,
the ideal generated by $X_1,\ldots,X_m$. 
The first resolution is minimal as  
$A/\rad (A)\cong kG\cong V(\ul)/\rad (V(\ul))$ as $A$-modules.
For minimality of the second resolution, note that the socle of $A$,
$\soc(A)$, is the
$k$-linear span of all $X_1^{\ell -1}\cdots X_m^{\ell -1}g_1^{b_1}\cdots
g_m^{b_m}$, where $0\leq b_i\leq \ell -1$.
We claim that in the notation of 
Lemma \ref{lem:free}, the socle of $V(\ul)$ has a basis in one-to-one
correspondence with the subset 
$$\{X_2^{\ell -1}\cdots X_{m}^{\ell -1} g_1^{b_1}\cdots
g_m^{b_m}\mid 0\leq b_i\leq \ell -1\}$$ of $\mathcal B$. 
Clearly $X_2,\ldots,X_m$ act trivially on all elements 
$X_2^{\ell -1}\cdots X_m^{\ell -1} g_1^{b_1}\cdots g_m^{b_m}\ot 1$
(in the notation of part (i) of this lemma).
We will check that $X_1$ also acts trivially:
\begin{eqnarray*}
  X_1X_2^{\ell -1}\cdots X_m^{\ell -1} g_1^{b_1}\cdots g_m^{b_m}\!
  \ot\! 1 \!\!\!&=&\!\!\! q^{-b_1} X_2^{\ell -1}\cdots X_m^{\ell -1} g_1^{b_1}\cdots
  g_m^{b_m} (\tau_{\ul}(t) \! -\! \sum_{i=2}^m \lambda_i X_i h_i)\!\ot\! 1\\
  \!\!\!& =&\!\!\! 0.
\end{eqnarray*}
It follows that $\soc(V(\ul))\cong kG \cong \soc(A)$ as $A$-modules.

(v) This follows immediately from (iv).
\end{proof}

We are also interested in simple $A$-modules.
The quotient $A/\rad(A)\cong kG$
is a commutative semisimple algebra. Thus the simple $A$-modules are all
one-dimensional, and correspond to the irreducible characters of $G$,
with $\Lambda$ acting trivially.

We shall use the notation 
$\underline{\Hom}$ to denote morphisms in the stable module category. 
In other words, 
$$\underline{\Hom}_A(M,N) = \Hom_A(M,N)/\PHom_A(M,N)$$
where $\PHom_A(M,N)$ is the set of all $A$-homomorphisms $f:M \to N$ which factor through a projective $A$-module. 
Recall that 
$$
 \Ext^n_A(M,N)\cong \underline{\Hom}_A(\Omega^n(M),N)\cong
  \underline{\Hom}_A(M,\Omega^{-n}(N)),
$$
where $\Omega^n$ (respectively, $\Omega^{-n}$) is the composition of
$n$ copies of $\Omega$ (respectively, $\Omega^{-1}$).
The following lemma will be needed in Section \ref{sec:ident}.

\begin{lem}\label{lem:a}
Let $S$ be a simple $A$-module. Then $\Ext^n_A(S,V(\underline{\lambda}))\neq 0$
for each $n$, $\ul$, and
the restriction map $\tau_{\ul}^*: \Ext^n_A
(S,V(\underline{\lambda})) \rightarrow \Ext^n_{k[t]/(t^{\ell})
}(S,V(\underline{\lambda}))$ is injective.
\end{lem}

\begin{proof}
Since  $A$ is free as a $\kt$-module by Lemma \ref{lem:free},
an $A$-injective resolution of $V(\ul)$ restricts to a
$\kt$-injective resolution.
It follows that $\Omega^{-n}_{\kt}(V(\ul))$ is isomorphic to
$\Omega^{-n}_A(V(\ul))$ in the stable module category, that is up
to projective direct summands. Thus
$$\Ext^n_{\kt}(S,V(\ul))\cong
  \underline{\Hom}_{\kt}(S,\Omega^{-n}_{\kt}(V(\ul))
 \cong \underline{\Hom}_{\kt}(S,\Omega^{-n}_A(V(\ul)).
$$

As $A$ is self-injective, $\Omega$ and $\Omega^{-1}$ are inverse operators
up to projective direct summands, so
by Lemma \ref{lem:vlambda}(v), $\Omega^{-n}_A(V(\ul))=V'(\ul)$ if $n$ is
odd, and $\Omega^{-n}_A(V(\ul))=V(\ul)$ if $n$ is even.
Assume without loss of generality that $\lambda_1 =1$.
Since $\soc(V(\ul))=kG X_2^{\ell -1}\cdots X_m^{\ell -1}\ot 1$ as a 
submodule of $V(\ul)\cong k\uparrow^A_{k\langle\tau_{\ul}(t)\rangle}$
(see the proof of Lemma \ref{lem:vlambda}(iv)), there is a unique (up to
scalar) nonzero $A$-homomorphism $f$ from $S$ to $V(\ul)$, sending $S$ to 
$ke_S X_2^{\ell -1}\cdots X_m^{\ell -1}\ot 1\subset \soc(V(\ul))$ 
where $e_S$ is the primitive central idempotent of $kG$ 
corresponding to $S$. This does not factor through a projective
$A$-module: If it did, it would factor through
$\displaystyle{A\stackrel{\cdot\tau_{\ul}(t)^{\ell -1}}
{\relbar\joinrel\relbar\joinrel\relbar\joinrel\longrightarrow}V(\ul)}$ 
since $A$ surjects onto $V(\ul)$.
The image of $S$ in $A$ must be contained in the socle of $A$, however 
the map $\cdot \tau_{\ul}(t)^{\ell -1}$ sends $\soc (A)$ to $0$.
Therefore this map represents an $A$-homomorphism from $S$ to $V(\ul)$ 
that is nonzero in $\underline{\Hom}_{A}(S,V(\ul))$.
A similar argument applies to $V'(\ul)$, proving that $\Ext^n_A(S,V(\ul))\neq 0$
for each $n$.

Next we show that the image of the map $f$ above, under restriction $\tau_{\ul}^*$,
remains nonzero in $\underline{\Hom}_{\kt}(S,V(\ul))$.
Again, if it does not, then $f: S\rightarrow V(\ul)$ factors as a
$\kt$-map through $\displaystyle{A\stackrel
{\cdot \tau_{\ul}(t)^{\ell -1}}{\relbar\joinrel\relbar\joinrel
\relbar\joinrel\longrightarrow} V(\ul)}$.
The image of $S$ in $A$ must be a one-dimensional
$k\langle \tau_{\ul}(t)\rangle$-submodule,
spanned by an element $a\in A$ for which $\tau_{\ul}(t) a =0$. 
Since $f$ sends a generator 
of $S$ to a non-zero element in $ke_S X_2^{\ell -1}\cdots X_m^{\ell -1}\ot 1 \subset \soc V(\ul)$, we get that
$a\tau_{\ul}(t)^{\ell -1}\in k^{\times}
e_S X_2^{\ell -1}\cdots X_m^{\ell -1} \tau_{\ul}(t)^{\ell -1}
=k^{\times} e_S X_1^{\ell -1}\cdots X_m^{\ell -1}$  under
the identification of $V(\ul)$ with $\Span_k ({\mathcal B} \tau_{\ul}(t)
^{\ell -1})\subset A$ in the notation of Lemma \ref{lem:free}. 
By Lemma \ref{compute} of the appendix, this cannot happen.
Hence, $f$ does not factor
through $\displaystyle{A\stackrel{\cdot \tau_{\ul}(t)^{\ell -1}}
{\relbar\joinrel\relbar\joinrel\relbar\joinrel\longrightarrow} V(\ul)}$.
A similar argument applies in odd degrees, involving $V'(\ul)$.
\end{proof}

\quad

Define an action of $G$
on projective space ${\mathbb{P}}^{m-1}$ by
\begin{equation}\label{gaction}
g_1^{a_1}\dots g_m^{a_m} \cdot [\lambda_1 : \dots : \lambda_m]
  = [q^{a_1}\lambda_1 : \dots : q^{a_m}\lambda_m].
\end{equation}

\begin{lem}\label{lem:projective}
Let $M$ be a finitely generated $A$-module. 
\begin{itemize}
\item[(i)] $\underline{\Hom}_{A}
(V(\ul),M)=0$ if, and only if, the restriction $M\!\downarrow
_{k\langle \tau_{\ul}(t)\rangle}$ is projective as a $k\langle
\tau_{\ul}(t)\rangle$-module. 
\item[(ii)] For each $g\in G$,
$M\!\downarrow_{k\langle\tau_{\ul}(t)\rangle}$ is
projective if, and only if, $M\!\downarrow
_{k\langle \tau_{g\cdot\ul}(t)\rangle }$ is projective.
\end{itemize}
\end{lem}

\begin{proof}
Lemma \ref{lem:vlambda}(i) together with the Eckmann-Shapiro Lemma 
implies the isomorphism
$$
  \underline{\Hom}_A(V(\ul),M)\cong \underline\Hom_A(k
  \uparrow_{k\langle\tau_{\ul}(t)\rangle}^A, M) \cong \underline{\Hom}_{k\langle
  \tau_{\ul}(t)\rangle}(k,M).
$$
This proves (i) 
since $\underline{\Hom}_{k\langle \tau_{\ul}(t)\rangle}(k,M) = 0$ if, and only if, 
$M\downarrow_{k\langle\tau_{\ul}(t)\rangle}$ is projective.
For (ii), note that $\tau_{g\cdot\ul}(t) = g\cdot \tau_{\ul}(t)$.
Since $g$ defines an inner automorphism of $A$, we now have 
$V(g\cdot\ul)\cong g\cdot V(\ul)\cong V(\ul)$.
Thus the statement follows from (i). 
\end{proof}

\section{Rank varieties}

In this section we define rank varieties for $A$-modules in the spirit of \cite{C1}.  The subalgebras 
${k\langle \tau_{\underline{\lambda}}(t)\rangle}$, 
defined in the text following (\ref{defn:taulambda}),
will play the role 
of cyclic shifted subgroups of $A$.

\begin{lem}\label{closed}
Let $M$ be a finitely generated $A$-module. The subset of projective space 
${\mathbb P}^{m-1}$,
consisting of all points $\underline{\lambda}$ such that $M\!\downarrow
_{k\langle \tau_{\underline{\lambda}}(t)\rangle}$ is not projective, is closed
in the Zariski topology.
\end{lem}

\begin{proof}
Let $n=\dim M$ and $S(\underline{\lambda})\in M_n(k)$ a matrix
representing the action of $\tau_{\underline{\lambda}}(t)$ on $M$.
Then $M$ is projective (equivalently, free) as a $k\langle
\tau_{\underline{\lambda}}(t)\rangle$-module if and only if the
Jordan form of $S(\underline{\lambda})$ has $n/\ell$ blocks of 
size $\ell$. That is $S(\underline{\lambda})$ has the maximal
possible rank for an $\ell$-nilpotent matrix, $n-n/\ell$.
The subset of $\mathbb P^{m-1}$,
$$
  \{\underline{\lambda}\in \mathbb P^{m-1}\mid \tau_{\underline{\lambda}}(t) \mbox{ 
does not have rank }n-n/\ell\},
$$
is described by the equations produced by the minors of $S(\underline
{\lambda})$ of size $(n-n/\ell)\times (n-n/\ell)$. All these
minors must be 0, and they give homogeneous  polynomial equations in the
coefficients $\lambda_i$ of $X_ih_i$. Thus this subset
is defined by a set of homogeneous  polynomials and is therefore closed.
\end{proof}

The action of $G$ by automorphisms on the polynomial algebra $k[X_1, \dots, X_m]$, defined by 
$g_i \cdot X_j = q^{\delta_{ij}}X_j$, is free and is easily seen to have the invariants
$$k[X_1, \dots, X_m]^G = k[X_1^{\ell}, \dots, X_m^{\ell}] \cong 
k[X_1, \dots, X_m].$$ 
Thus, $\mathbb A^m/G  = \Spec k[X_1, \dots, X_m]^G \cong \mathbb A^m$,
where $\mathbb A^m$ is the affine space $k^m$ (see for example \cite[I.5.5(6)]{Jan} for the first equality).
Since the action of $G$ commutes with the standard action of $k^*$ on $\mathbb A^m$, and the induced action on $\mathbb P^{m-1}  = \mathbb A^m/k^*$ 
is  the action defined as in (\ref{gaction}), we have
$\mathbb P^{m-1}/G \cong \mathbb P^{m-1}$.
Furthermore, 
Lemma \ref{lem:projective}(ii) implies that the set 
$
\{\ul\in{\mathbb{P}}^{m-1}\mid M\downarrow
   _{k\langle\tau_{\ul}(t)\rangle} \mbox{ is not projective}\}
$
is stable under the action of $G$. 
Thus, we can make the following definition.

\begin{defn}
\label{defn:rank}
The {\em rank variety} of an $A$-module $M$ is
$$
  V_A^r(M)=\{\ul\in{\mathbb{P}}^{m-1}\mid M\downarrow
   _{k\langle\tau_{\ul}(t)\rangle} \mbox{ is not projective}\}/G.
$$\end{defn} 

\noindent
We will sometimes abuse notation and write $\ul\in V^r_A(M)$ when we mean
that $\ul$ is a representative of a $G$-orbit in $V^r_A(M)$.
Note that Lemma \ref{closed} ensures $V^r_A(M)$ is a projective
variety for any finitely generated  $A$-module $M$. The following properties of these varieties are immediate.

\begin{prop}
\label{prop:rank}
Let $M, N, M_1, M_2, M_3$ be $A$-modules.

\begin{itemize}
\item[(i)] $V^r_A(k) = \mathbb P^{m-1}/G \cong \mathbb P^{m-1}$.
\item[(ii)] $V^r_A(M\oplus N)=V^r_A(M)\cup V^r_A(N)$.
\item[(iii)] $V^r_A(\Omega^i (M))=V^r_A(M)$ for all $i$.
\item[(iv)] If $0\rightarrow M_1\rightarrow M_2\rightarrow M_3\rightarrow 0$
is a short exact sequence of $A$-modules, then 
$V^r_A(M_i)\subset V^r_A(M_j)\cup V^r_A(M_k)$ for any $\{i,j,k\}=
\{1,2,3\}$.
\end{itemize}

\end{prop}

We will denote $V_A^r(k)$ by $V_A^r$.

\quad

The rank variety characterizes projectivity of modules by the following
lemma, a version of Dade's Lemma for finite group representations
\cite{Dade}.
We thank K.~Erdmann and D.~Benson for suggesting to us that the proof of a 
generalization of Dade's Lemma in \cite{BCR} should apply almost verbatim in our setting. 
For completeness, we give our
adaptation of the proof in \cite{BCR} here (cf.\ \cite[Thm.\ 2.6]{BEH}).

\begin{thm} \label{thm:dade}
Let $M$ be a finitely generated $A$-module.
Then $V^r_A(M)=\emptyset$ if, and only if, $M$ is a projective $A$-module.
\end{thm}

\begin{proof}
If $M$ is projective, then $M\downarrow_{k\langle\tau_{\ul}(t)\rangle}$
is projective for all $\ul$ by Lemma \ref{lem:free}, so
$V^r_A(M)=\emptyset$.

For the converse, we argue by induction on $m$. 
Let $Y_i = X_ih_i$ ($i=1,\ldots,m$), where $h_i$ is defined in the text
preceding (\ref{eqn:q-commute}).
Let $\Lambda_m'=k\langle Y_1,\ldots,Y_m\rangle$ and note that $M$ is
projective if and only if $M\!\downarrow_{\Lambda_m'}$ is projective:
We may write $A\cong \Lambda_m'\rtimes G$.
If $M\!\downarrow_{\Lambda_m'}$ is projective, any surjective $A$-map from
another $A$-module $N$ onto $M$ splits on restriction to $\Lambda_m'$.
The splitting map may be averaged by applying $\frac{1}{|G|}\sum_{g\in G}
g$ to obtain an $A$-map, as the characteristic of $k$ does not divide
$|G|$.  
If $m=1$, this immediately implies that $M$ is projective 
if, and only if, $V^r_A(M)=\emptyset$.

Let $m=2$ and assume $V^r_A(M)=\emptyset$ but that $M$ is not projective.  
We will show that $M$ is forced to be $0$.
The Jacobson radical of $A$ is $J = (Y_1, Y_2)$, the ideal
generated by $Y_1,Y_2$. Let 
$$
  N = \{ u \in M | Ju \subset J^{\ell-1}M \}.
$$Let $Y = \lambda_1Y_1+\lambda_2Y_2=\tau_{\ul}(t)$.  
We will first show that the map induced by $Y$:
$$
N/J^{\ell-1}M \stackrel{\cdot Y}{\relbar\joinrel\rightarrow} 
J^{\ell-1}M/J^{\ell} M
$$
is an isomorphism for any pair $(\lambda_1, \lambda_2)$ where $\lambda_2 
\not = 0$.  
We will need the observation that $YJ^{\ell -1}=J^{\ell}$, which follows
from $Y_1^{\ell -i}Y_2^i\in YJ^{\ell -1}$ for $i\in\{1,\ldots,\ell -1\}$
as may be proven by induction on $i$: 
If $i=1$, then $YY_1^{\ell -1}=\lambda_2 q^{\ell -1}
Y_1^{\ell -1}Y_2$, so $Y_1^{\ell -1}Y_2\in YJ^{\ell -1}$. If $i\geq 2$,
then $YY_1^{\ell -i}Y_2^{i-1}=\lambda_1Y_1^{\ell -i+1} Y_2^{i-1}
+\lambda_2 q^{\ell -i} Y_1^{\ell -i}Y_2^i$, so
$Y_1^{\ell -i}Y_2^i\in YJ^{\ell -1}$ by induction.

\vspace{0.05in}

Injectivity of $\cdot Y$: Let $v \in N$. 
Suppose $Yv \in J^{\ell}M=YJ^{\ell -1}M$.  
Then there exists $u \in J^{\ell-1}M$ such that 
$Yv = Yu$. Therefore, $Y(v-u) = 0$. Since $M\!\downarrow_{k\langle Y \rangle}$
is projective, we have 
$v-u = Y^{\ell-1} u^\prime$ for some $u^\prime \in M$. 
Hence, $v = u + Y^{\ell-1} u^\prime \in J^{\ell-1}M$. In other words, 
$\bar v = 0 \in N/J^{\ell-1}M$. We conclude that $\cdot Y$ is injective.

\vspace{0.05in}
Surjectivity of $\cdot Y$: 
We may assume that $M$ does not have projective summands.  
This implies that $\soc(k[Y_1, Y_2]/(Y_1^{\ell}, Y_2^{\ell}))M = 
Y_1^{\ell-1}Y_2^{\ell-1}M = 0$.
Now the relations on $Y_1,Y_2$ easily imply that 
\begin{equation}
\label{socle}
Y^{\ell -1}Y^{\ell -1}_1 M =  Y_1^{\ell -1}Y_2^{\ell-1}M = 0.
\end{equation}
To show surjectivity we need to show that $YN = J^{\ell-1}M$. 
Observe that $J^{\ell-1} = kY_1^{\ell-1} + YJ^{\ell-2}$. Thus,   
$J^{\ell-1}M = Y_1^{\ell-1}M + YJ^{\ell-2}M$. 
The definition of $N$ immediately implies 
that $J^{\ell-2}M \subset N$. Therefore
$YJ^{\ell-2}M \subset YN$. 
Hence, to show the inclusion $Y_1^{\ell-1}M + YJ^{\ell-2}M \subset YN$, it 
suffices to show that $Y_1^{\ell-1}M \subset YN$.

Take any element of the form $Y^{\ell-1}_1 u$, $u\in M$.  
By (\ref{socle}), $Y^{\ell-1}Y^{\ell-1}_1 u =0$.  Since 
$M\downarrow_{k\langle Y \rangle}$ is projective, there is an
element $u'\in M$ such that
\begin{equation}
\label{first}
Y_1^{\ell -1}u = Yu^\prime.
\end{equation}  
Multiplying both sides by $Y_1$, we get $Y_1Yu^\prime = 0$. Thus, 
$Y_1(\lambda_1Y_1+\lambda_2Y_2)u^\prime = 0.$
Using the relation $Y_2Y_1=qY_1Y_2$, we get 
$
 (\lambda_1Y_1 + q^{-1}\lambda_2 Y_2)Y_1 u^\prime = 0.
$
Applying our projectivity 
hypothesis to the restriction of $M$ to $k\langle\lambda_1Y_1+
q^{-1}\lambda_2Y_2
\rangle$, there is a $u^{\prime\prime}\in M$ for which
\begin{equation}
\label{second}
Y_1 u^\prime = (\lambda_1Y_1+q\lambda_2Y_2)^{\ell -1}u^{\prime\prime}.
\end{equation}
Combining (\ref{first}) and (\ref{second}), we get 
$$
Yu^\prime \in J^{\ell -1}M, \quad 
Y_1 u^\prime \in J^{\ell -1} M.$$
Since $\lambda_2 \not = 0$, $Y$ and $Y_1$ generate $J$. Thus,
$Ju^\prime \subset J^{\ell -1}M$, so that $u^\prime\in N$ by definition.
By (\ref{first}) it follows that $Y_1^{\ell -1}M\subset YN$, as was needed.

\vspace{0.1in} 

Thus, we obtain that for any non-zero pair $(\lambda_1, \lambda_2)$, 
the map 
$$
\cdot(\lambda_1Y_1+\lambda_2Y_2): N/J^{\ell -1}M \to J^{\ell -1}M/J^{\ell} M
$$
is an isomorphism.   Assume $J^{\ell -1}M/J^{\ell} M \not = 0$.
Use the isomorphism $\cdot Y_2$ to
identify $N/J^{\ell -1}M$ with $J^{\ell -1}M/J^{\ell} M$,
so that we may consider the maps $\cdot (\lambda_1Y_1+\lambda_2Y_2)$ to be
endomorphisms of $J^{\ell -1}M/J^{\ell}M$.
Taking $\lambda_2=1$, the determinant of the operator $\cdot
(\lambda_1Y_1+Y_2)$ is a polynomial in $\lambda_1$ and thus there is a
value of $\lambda_1$ for which the operator is not invertible, a
contradiction.
Thus $J^{\ell -1}M/J^{\ell} M = 0$, and Nakayama's Lemma implies that $J^{\ell 
-1} M = 0$. 
In particular, $Y_2^{\ell -1} M = 0$. 
At the same time, $M$ restricted to $k\langle Y_2\rangle$ is projective.  
Therefore $M = 0$, completing the proof in the case $m=2$.  
\vspace{0.1in}

Now suppose $m \geq 3$ and $V^r_A(M)=\emptyset$. 
We will show that $M\!\downarrow_{\Lambda_m'}$ is projective, where
$\Lambda_m'=k\langle Y_1,\ldots,Y_m\rangle$.
As noted at the beginning of the proof, this will imply $M$
is projective.
  
We have a short exact sequence of algebras 
(in the sense of \cite[XVI \S 6]{CE}):
$$
k\langle Y_m \rangle \to \Lambda^\prime_m \to \Lambda^\prime_{m-1}.
$$
Therefore there is a spectral sequence
$$
\coh_p(\Lambda^\prime_{m-1}, \coh_q(k\langle Y_m \rangle,M)) \Rightarrow
\coh_{p+q}(\Lambda^\prime_m,M).
$$
By our assumption, $M\downarrow_{k\langle Y_m \rangle}$ is projective. 
Thus the spectral sequence collapses at $E^2$ and we get an isomorphism
$$
\coh_p(\Lambda^\prime_{m-1}, M/Y_m M) \cong    \coh_p(\Lambda^\prime_m,M).
$$
Therefore, to finish the proof it suffices to check that $M/Y_m M$ is 
projective as a $\Lambda^\prime_{m-1}$-module.
 
Write $\ul^\prime = [\lambda_1:\cdots : \lambda_{m-1}]\in{{\mathbb P}}^{m-2}$,
and let $\tau^\prime(\ula^\prime) = \lambda_1Y_1 + \dots +
\lambda_{m-1}Y_{m-1}$. 
Consider the subalgebra $B\subset \Lambda^\prime_m$ generated by 
$\tau^\prime(\ula^\prime)$ and $Y_m$. Any element of $B$ of the form 
$\mu_1 \tau^\prime(\ula^\prime) + \mu_2 Y_m$ is of the form 
$\tau(\ula)$ for the algebra $\Lambda^\prime_m$. Thus, 
$M\!\downarrow_{k\langle\mu_1\tau^\prime(\ula^\prime) + \mu_2Y_m\rangle}$ 
is projective for 
any pair $(\mu_1, \mu_2)$. 
 Since $\tau^\prime(\ula^\prime)$ and $Y_m$ $q$-commute, the argument for 
$m=2$ applies to $B$. 
 We conclude that $M$ is projective as a $B$-module. 
Therefore, $M/Y_mM$ is  projective as $B/Y_m B$-module.
In other words, 
 $M/Y_m M\downarrow_{k\langle\tau^\prime(\ula^\prime) \rangle}$ is
projective.  
By the induction hypothesis, $M/Y_mM$ is projective as a 
$\Lambda_{m-1}^\prime$-module.   
 Thus, $\coh_{p}(\Lambda^\prime_m,M) = 0$ for all $p>0$.  
Since the trivial module is the only
simple module for the local algebra $\Lambda^\prime_m$, we conclude that
 $M$ is projective as a $\Lambda^\prime_m$-module. Therefore, $M$ is
projective as an $A$-module.
  \end{proof}

\section{Support varieties}

In this section we introduce cohomological support varieties for $A$-modules. 
The proofs of their properties are standard and will be omitted when they are identical to those existing in the literature.

The cohomology of the quantum elementary abelian group $A$ is
\begin{equation}\label{eqn:cohA}
 \coh^*(A,k) = \Ext^*_A(k,k)\cong k[y_1,\ldots,y_m]
\end{equation}
as a graded algebra, where $\deg(y_i)=2$. 
Indeed, let $A_1 = (k[t]/(t^{\ell}))\!\rtimes\!\Zl$, the algebra $A$ in the case $m=1$. 
The periodic 
$k[t]/(t^{\ell})$-free resolution of $k$,
\begin{equation}\label{eqn:tres}
  \cdots \stackrel{\cdot t^{\ell -1}}{\relbar\joinrel
  \longrightarrow} k[t]/(t^{\ell})\stackrel{\cdot t}{
  \longrightarrow}k[t]/(t^{\ell}) 
  \stackrel{\cdot t^{\ell -1}}{\relbar\joinrel
  \longrightarrow} k[t]/(t^{\ell})\stackrel{\cdot t}{
  \longrightarrow}k[t]/(t^{\ell})\stackrel{\varepsilon}{\longrightarrow}
  k \rightarrow 0,
\end{equation}
becomes an $A_1$-projective resolution by giving $k[t]/(t^{\ell})$ the
standard $\Zl$-action $g\cdot t^i=q^it^i$ in even degrees and the shifted 
$\Zl$-action $g\cdot t^i=q^{i+1}t^i$  in odd degrees,
where $g$ is a generator of $\Zl$. 
The resolution yields $\coh^*(A_1,k) \cong k[y]$. 
The general case is obtained by applying the K\"unneth formula.

Recall that $\coh^*(A,M)=\Ext^*_A(k,M)$ is an $\coh^*(A,k)$-module under Yoneda composition, for any $A$-module $M$.
To proceed with our geometric constructions, we will need to establish finite generation of  
$\coh^*(A,M)$  over $\coh^*(A,k)$ whenever $M$ is finitely generated.

\begin{lem}
\label{finite}
Let $M$ be a finitely generated $A$-module.  
Then $\coh^*(A,M)$ is finitely generated as an $\coh^*(A,k)$-module.
\end{lem}

\begin{proof}
Since  $M$ is finitely generated, 
we can argue by induction on the length of its composition series.  
Hence, it suffices to prove the lemma for simple $A$-modules.  
Let $S$ be such a module.
The spectral sequence in cohomology arising from the sequence of augmented 
algebras (see, for example, \cite[XVI \S6]{CE}) $ \Lambda \to A \to kG $
yields the isomorphism 
$\coh^*(A,S) = \coh^*(\Lambda, S)^G$.
Let $R = \coh^*(\Lambda,k)$, a finitely generated $G$-algebra.  
Since any simple $A$-module becomes trivial when restricted to $\Lambda$, 
we conclude that $\coh^*(\Lambda, S)$ is a rank 1 free 
$R$-module with a compatible action of $G$.  
Since $\coh^*(A,k) = \coh^*(\Lambda, k)^G = R^G$,  
it remains to see that $\coh^*(\Lambda, S)^G$ is finitely generated over $R^G$.
As the characteristic of $k$ does not divide the order of $G$, this is a consequence of
the Noether Theorem stating that the Noetherian $k$-algebra $R$ is finitely generated
over $R^G$ (see, for example, \cite[1.3.1]{Ben2}). 
\end{proof}

\begin{rem}
\label{simple}
We can compute $\coh^*(A,S)$ explicitly when $S$ is simple,
yielding more insight in
our special case.  First let $m=1$. Let $S_i$ be the (one-dimensional) simple $A_1$-module
on which $t$ acts as multiplication by 0 and the generator $g$ of $\Zl$
acts as multiplication by $q^i$.
Using the resolution (\ref{eqn:tres}), we get 
$$
  \coh^n(A_1,S_i)=\left\{
  \begin{array}{ll}
   0, & \text{for all }n\text{ if } i\not\in \{0,1\}\\
  0, & \mbox{if }n\mbox{ is even and }i=1, \mbox{ or if }n\mbox{ is odd and }i=0 \\
  k, & \mbox{if }n\mbox{ is odd and }i=1, \mbox{ or if }n\mbox{ is even and }i=0.\end{array}\right.
$$
Since the action of the generator $y$ of $k[y] \cong \coh^*(A_1,k)$ induces a periodicity 
isomorphism, we get that $\coh^*(A_1, S_i)$ is a rank $1$ free $\coh^*(A_1,k)$-module in
case $i\in \{0,1\}$.

Now if $m\geq 2$, a (one-dimensional) simple $A$-module may be written
$S=S_{\chi}$ for some $\chi:G\rightarrow k^{\times}$, where
any $g\in G$ acts as multiplication by $\chi(g)$.
It may be factored as $S_{\chi}\cong S_{\chi_1}\otimes \cdots \otimes S_{\chi_m}$
where $g_j$ acts trivially on $S_{\chi_i}$ if $i\neq j$.
We may similarly factor $A\cong A_1\otimes\cdots\otimes A_m$ where 
$A_i= (k[X_i]/(X_i^{\ell}))\rtimes\langle g_i\rangle$.
Apply the K\"unneth Theorem to obtain
$$ 
  \coh^*(A,S_{\chi})\cong \coh^*(A_1,S_{\chi_1})\otimes\cdots\otimes
  \coh^*(A_m,S_{\chi_m}),
$$
and apply the case $m=1$ to each factor.
If $k = \mathbb C$, this result and Lemma \ref{finite} follow  from 
\cite[Thm.\ 2.5]{GK}, since 
 in this case $A$ is isomorphic to the quantized restricted Borel subalgebra of $u_q(sl_2 ^{\times m})$.
\end{rem}

\vspace{.1in}

Let $M,N$ be left $A$-modules. Let $I(M,N)$ be the annihilator of
$\Ext^*_A(M,N)$ under the action of 
$\coh^*(A,k)$ by cup product, equivalent to $-\ot _k N$ followed
by Yoneda composition (see \cite[I.\ Prop.\ 3.2.1]{Ben}).
Since $I(M,N)$ is a homogeneous ideal, it defines a projective subvariety of 
$\Proj \coh^*(A,k) \cong \mathbb P^{m-1}$. 
We define
\begin{equation}\label{eqn:suppMN}
V_A^c(M,N) = \Proj  \coh^*(A,k)/I(M,N),
\end{equation}
the set of all homogeneous prime ideals of dimension $1$ (that is,
those ideals not contained in  any homogeneous prime ideal other than 
$\coh^{*>0}(A,k)$) which contain $I(M,N)$. 
If $M=N$, we write simply $I(M)=I(M,M)$ and $V_A^c(M)=V_A^c(M,M)$. 
Note that 
$$V_A^c(k)= \Proj  \coh^*(A,k) \cong {\mathbb{P}}^{m-1}.$$
We will denote $V^c_A(k)$ by $V^c_A$.

\begin{prop}
 \label{prop:prop} The following properties hold for all finitely generated $A$-modules $M,N$.
\begin{itemize}
\item[(i)] $V_A^c(M \oplus N) = V_A^c(M) \cup V_A^c(N)$.
  \item[(ii)] $\va(M) = \va(\Omega (M))$.
\item[(iii)] $\va(M,N) \subset \va(M) \cap \va(N)$.
\item[(iv)] $V^c_A(M)=\cup_S V^c_A(S,M)$, where $S$ runs over all simple $A$-modules.
\item[(v)] $V^c_A(M) = \emptyset$ if, and only if, $M$ is projective.
\item[(vi)] $V_A^c(M \otimes N) \subset V_A^c(M) \cap V_A^c(N)$.

\end{itemize}
\end{prop}

\begin{proof}
Arguments from \cite[II.\ \S5.7]{Ben} apply verbatim to prove (i)--(iv).

(v) If $M$ is projective, then $\Ext^{n}_A(M,M) = 0$ for all $n>0$. 
Thus, $V^c_A(M) = \emptyset$.

Assume $V^c_A(M) = \emptyset$.  
By (iv), we get $V^c_A(S,M) = \emptyset$ for every simple $A$-module $S$. 
Note that $\Ext_A^*(S,M) \cong \coh^*(A, S^\#\otimes M)$ (see \cite[I.\ \S3.1]{Ben}
or the appendix), and Lemma \ref{finite} thus implies that this cohomology
is finitely generated as an $\coh^*(A,k)$-module.
So there exists $n_0$ such that $\Ext^n_A(S,M)=0$ for any $S$ and any
$n>n_0$.
Hence, the minimal injective resolution of $M$ is finite.   
Since injective $A$-modules are projective, 
any finite resolution splits. Hence, $M$ is injective, and therefore  projective.  

\vspace{0.1in}
\noindent
(vi) The action of $\coh^*(A,k)$ on $\Ext^*_A(M\otimes N,M\otimes N)$
factors through its action on $\Ext^*_A(M,M)$: 
We may first apply $ - \otimes M$ to an $n$-extension of $k$ by $k$,
and then apply $ - \otimes N$. Thus $I(M,M)\subset I(M\otimes N,M\otimes N)$,
implying that $V^c_A(M\otimes N)\subset V^c_A(M)$. 
On the other hand, since $M$ is finitely generated, 
we have the adjunction isomorphism 
$\Ext^*_A(M\otimes N,M\otimes N) \cong \Ext^*_A(M^\#\otimes M \otimes N, N)$ 
(see the appendix). 
Thus, $\va(M\otimes N) = \va(M^\#\otimes M \otimes N, N)$. The latter 
is contained in $\va(N)$ by (iii).
\end{proof}

Following the original construction by Carlson \cite{C2} 
for finite groups, we introduce modules $L_\zeta$. 
Let $\zeta \in \coh^{n}(A,k) \cong \underline \Hom(\Omega^{n}(k),k)$. 
Then $L_\zeta$ is defined to be the kernel of the corresponding map
$\zeta: \Omega^{n}(k) \to k$. In other words, we have a short exact sequence
\begin{equation}\label{exact}
0 \to L_\zeta \to \Omega^{n}(k) \to k \to 0.
\end{equation}
Let $\langle \zeta \rangle \subset \mathbb P^{m-1}$ be the projective hypersurface defined by $\zeta$, that is
the set of all homogeneous prime ideals of dimension 1 which contain 
$\zeta$.

The following result is an adaptation to our situation of \cite[3.3]{EHSST}.
\begin{prop}
$V^c_A(M \otimes L_\zeta) = \langle\zeta\rangle \cap V_A^c(M)$.
\end{prop}

\begin{proof}  
Observe that for any $M, N$, a prime ideal $\wp$ belongs to $\va(M,N)$ if and only if $I(M,N) \subset \wp$ if and only 
if $\Ext_A^*(M,N)_\wp \not = 0$. 

We first show $ \langle\zeta\rangle \cap V_A^c(M) \subset V^c_A(M \otimes L_\zeta)$.
By Proposition \ref{prop:prop}(iv), 
$V_A^c(M \otimes L_\zeta) = \bigcup \va(S, M \otimes L_\zeta)$ and $\va(M) = \bigcup \va(S,M)$,
where $S$ runs through all simple $A$-modules,
so it suffices to show that
$$\langle\zeta\rangle \cap V_A^c(S,M) \subset \va(S, M \otimes L_\zeta)$$ for any $S$. 
Let $\wp$ \, be a homogeneous prime ideal in $\langle\zeta\rangle \cap V_A^c(S,M)$, that is
$\wp$ \, contains the ideal generated by $I(S,M)$ and $\zeta$. 
We want to show that $\wp \in V^c_A(S,M \otimes L_\zeta)$, that is $I(S, M \otimes L_\zeta) \subset \wp$.

 Suppose $I(S, M \otimes L_\zeta) \not \subset \wp$.  
This implies that $\Ext_A^*(S, M \otimes L_\zeta)_{\wp} = 0$. 
Tensoring the short exact sequence $0 \to L_\zeta \to \Omega^{n}( k) \to k \to 0$ with $M$ and applying $\Ext_A^*(S, -)$, 
we get a long exact sequence
$$
 \dots\! \to\! \Ext_A^i(S, M \otimes L_\zeta)\! \to\! \Ext_A^{i-n}(S,M)\! 
\stackrel{\zeta}{\to}\! \Ext_A^i(S,M)\! 
 \stackrel{\delta}{\to}\! \Ext_A^{i+1}(S, M \otimes L_\zeta)\! \to\! \dots 
 $$
Let $z \in \Ext^i(S,M)$.  Then $\delta(z) \in  \Ext_A^{i+1}(S, M \otimes L_\zeta)$. Since $\Ext_A^{i+1}(S, M \otimes L_\zeta)_\wp = 0$, 
there exists a homogeneous element $a \not \in \wp$ such that $
\delta(az) = a\delta(z)= 0$.
Moreover, we can assume that the cohomological degree of $a$  is sufficiently large, so that
$\deg(a)i > n$. 
 The long exact sequence implies that $az = \zeta y$ for $y \in \Ext_A^{\deg(a)i-n}(S,M)$.  We conclude that 
 $$\Ext_A^*(S,M)_\wp = \zeta \Ext_A^*(S,M)_\wp .$$ 
Since $\zeta \in \wp$, and $\Ext_A^*(S,M) = \Ext_A^*(k, S^\#\otimes M)$ is finitely generated over $\coh^*(A,k)$ by Lemma \ref{finite}, 
Nakayama's Lemma implies that 
$\Ext_A^*(S,M)_\wp = 0$. But this contradicts the assumption $I(S,M) \subset \wp$.  We conclude that  $I(S, M \otimes L_\zeta) 
 \subset \wp$, 
and hence $\langle\zeta\rangle \cap V_A^c(S,M) \subset \va(S, M \otimes L_\zeta)$.

\vspace{0.1in}
To prove the opposite inclusion $V^c_A(M \otimes L_\zeta) \subset \langle\zeta\rangle \cap V_A^c(M) $, it suffices, 
by Proposition \ref{prop:prop}(vi), to show that $\va(L_\zeta) \subset  \langle\zeta\rangle$. Using \ref{prop:prop}(iv) again, 
we reduce to showing the inclusion $\va(S, L_\zeta) \subset \langle\zeta\rangle$ for any  simple $A$-module $S$.  
Thus we need to show that $\Ext^*_A(S,L_\zeta)_\wp \not = 0$ for a prime ideal $\wp\subset \coh^*(A,k)$ implies $\zeta \in \wp$. 
We will prove the converse. Suppose $\zeta \not \in \wp$.  Then multiplication by $\zeta$ induces an isomorphism on 
$\Ext^*(S,k)_\wp$ since it is invertible in $\coh^*(A,k)_\wp$. Since localization is exact, 
the short exact sequence (\ref{exact})  implies that $\Ext^*(S,L_\zeta)_\wp$ is the kernel of 
the isomorphism $\zeta: \Ext^*(S,k)_\wp \to \Ext^{*+n}(S,k)_\wp$. Thus, $\Ext^*(S,L_\zeta)_\wp = 0$.  
  \end{proof}

Induction yields the following corollary.

\begin{cor}
\label{cor:realization}
$V^c_A(M \otimes L_{\zeta_1}\otimes \dots \otimes L_{\zeta_i}) = \langle\zeta_1, \dots ,\zeta_i\rangle \cap V_A^c(M)$.
\end{cor}


\section{Identification of varieties}\label{sec:ident}

In this section we will establish an analogue of the Avrunin-Scott Theorem, 
identifying the cohomological variety with the rank variety of  a module.  
For $\ul \in \mathbb P^{m-1}$ we denote by 
$\tau_{\ul} ^*: \coh^*(A,k) \to \coh^{ev}(k[t]/(t^{\ell}), k) $ 
the map induced on cohomology by 
$\tau_{\ul} :k[t]/(t^{\ell}) \hookrightarrow A$ as defined in
(\ref{defn:taulambda}).

Recall the algebra $\widetilde{A}=k[X_1,\ldots,X_m]\rtimes G$ 
defined at the beginning of \S 2.
We have a short exact sequence of augmented algebras
\begin{equation}
\label{power}
        k[X_1, \dots, X_m] \stackrel{X_i \mapsto X_i^{\ell}}{\relbar\joinrel
\relbar\joinrel\relbar\joinrel\longrightarrow} \widetilde{A} \to  A
\end{equation}
which induces a spectral sequence (see \cite[XVI \S6]{CE})
\begin{equation}
\label{spseq}
\coh^*(A, \coh^*(k[X_1, \dots, X_m],k)) \Rightarrow \coh^*(\widetilde{A},k).
\end{equation}

\begin{lem}
\label{lem:spseq} 
The transgression map of the $E_2$ page of the spectral sequence (\ref{spseq})
induces an isomorphism $\coh^1(k[X_1, \dots, X_m],k) \cong \coh^2(A,k)$. 
\end{lem}

\begin{proof} 
We first show that the action of $A$ on $\coh^*(k[X_1,\ldots,X_m],k)$ is trivial. Since this action comes from 
tensoring the action of each factor
$(k[X_i]/(X^{\ell}_i)) \rtimes \Zl$ of $A$ on the corresponding $\coh^*(k[X_i],k)$, it suffices to check this for $m=1$. 
Let $\widetilde{A}_1 = k[X] \rtimes \Zl$, $A_1 = (k[X]/(X^{\ell})) \rtimes \Zl$.  Let $\mathcal L =k[X^\ell] \subset \widetilde A_1$ denote 
the image of $k[X]$ in $\widetilde A_1$. Since $\coh^i(k[X],k) = 0$ for $i>1$ and $\coh^0(k[X],k) = k$, the trivial module,
we need only check that the action of $A_1$ on $\coh^1(k[X],k)$ is trivial. 

The Koszul resolution $K_*: 0 \to k[X] \stackrel{\cdot X}{\to} k[X]
\stackrel{\varepsilon}{\to}k \to 0$ of $k$ as a $k[X]$-module
becomes a resolution as a $k[X]\rtimes\Zl$-module under the action $g \circ X^i = q^iX^i$ in degree $0$ and 
$g \circ X^i = q^{i+1}X^i$ in degree $1$.  The spectral sequence (\ref{spseq}) 
can be obtained as a 
Grothendieck spectral sequence associated to the composition of functors 
$\Hom_{\mathcal L}(-,k)$ and $\Hom_{A_1}(k,-)$. Hence, 
 the action of $A_1$ on $\coh^*(\mathcal L,k) = \Ext^*_\mathcal L(k,k)$ 
is induced by the action of $A_1$ on the complex 
$$\Hom_\mathcal L(K_*, k): \ \ 0 \leftarrow \Hom_k(k[X]/(X^\ell),k) \stackrel{\delta}{\leftarrow} \Hom_k(k[X]/(X^\ell),k) \leftarrow 0$$
which, in turn, comes from the diagonal
action of $\widetilde A_1$ on $\Hom_k(K_*,k)$ (see the appendix for the explicit formula).
Computation yields that $H^1(\mathcal L,k)$ is generated by 
the cocycle  $f:k[X]/(X^\ell) \to~k$ specified by the condition
$f(X^{\ell -1}) = 1, f(X^i) = 0$ for $i \not = \ell -1$.  We further compute 
that $ X \circ f \in \Ima \delta$ and $ (g \circ f)(X^{\ell-1}) = f(g^{-1}\circ X^{\ell-1}) = f(g^{\ell-1}\circ X^{\ell-1}) = 
f(q^\ell X^{\ell-1})= f(X^{\ell-1})$. Thus, the action of $A_1$ on $H^1(\mathcal L,k)$ is trivial.
\sloppy
{

}

Hence, $A$ acts trivially on $\coh^q(k[X_1, \dots, X_m],k)$. Therefore, 
$$E_2^{0,q} = \coh^q( k[X_1,\ldots,X_m], k)^A = \coh^q( k[X_1,\ldots,X_m], k).$$ 

Observe that $\coh^1(\widetilde A,k) = 0$: We have
$\coh^1(\widetilde A,k) \cong \bigoplus\limits_{m} 
\coh^1(k[t]\rtimes \Zl,k)$ by the K\"unneth Theorem. 
Direct
computation with the Koszul resolution $K_*$ of $k$ 
shows that $\coh^1(k[t]\rtimes \Zl , k)=0$.
Thus $\coh^1(\widetilde A,k) = 0$ as required.

It follows that $E_3^{0,1} = E_\infty^{0,1}=0$, and so $\Ker d_2^{0,1} = 0$. Thus 
$$d_2^{0,1}: \coh^1(k[X_1, \dots, X_m], k) \to \coh^2(A,k)$$ is injective.  
Since $\dim_k \coh^2(A,k) = m = \dim_k \coh^1(k[X_1, \dots, X_m],k)$, 
we conclude that $d_2^{0,1}$ is an isomorphism.
\end{proof}

The next lemma establishes that the map $\tau_{\ul}^*$ is ``essentially surjective" and is invariant under the $G$-action on $\mathbb P^{m-1}$.
Let $I = (X_1, X_2, \dots, X_m)$ be the ideal of  $A$ generated by the $X_i$'s, 
  and let $z_1, \dots, z_m$  be the basis of $(I/I^2)^\#$ dual to 
$X_1, \dots, X_m$.
As before, we denote by $y_i$ the generators of $\coh^*(A,k)$.  

\begin{lem}
\label{lem:well}
(i)  For any $\ul \in \mathbb P^{m-1}$, $\tau_{\ul}^*$ is surjective onto $\coh^{ev}(k[t]/(t^{\ell}), k) \cong k[y]$.
 
(ii) For  any $g \in G$,
$\tau_{\ul}^* = \tau_{g \cdot \ul}^*$. 
\end{lem}

\begin{proof}
(i) 
Let $\ula = [\lambda_1: \lambda_2: \cdots :\lambda_m] \in \mathbb P^{m-1}$. 
Lemma \ref{lem:tau} implies that the following diagram is commutative, 
where both rows are exact sequences of augmented algebras:
\begin{equation}
\label{lambda}
\begin{xy}*!C\xybox{
\xymatrix{
k[X_1, \dots, X_m] \ar[rrr]^{X_i \mapsto X_i^{\ell} \quad} &&& \widetilde A \ar[rrr] &&& A \\
k[t] \ar[rrr]^{t \mapsto t^{\ell}} \ar[u]_{t \mapsto \lambda_1^{\ell} X_1 + \dots + \lambda_m^{\ell} X_m}  
&&& k[t] \ar[rrr] \ar[u]_{t\mapsto \lambda_1 X_1h_1 + \dots + \lambda_m X_m h_m} &&& k[t]/(t^{\ell})  \ar[u]_{t \to \tau_{\ul}(t)}
}}
\end{xy}
\end{equation}

The rows of (\ref{lambda}) induce compatible spectral sequences where the edge homomorphisms $d_2^{0,1}$ are  
isomorphisms by Lemma \ref{lem:spseq}.
Thus, we get another commutative diagram where vertical maps are restrictions.
\begin{equation}
\label{lambdadual}
\begin{xy}*!C\xybox{
\xymatrix{
(I/I^2)^\# \ar@{=}[r] \ar[d] & \coh^1(k[X_1, \dots, X_m],k) \ar[d] \ar[r]^>>>>>{\sim}& \coh^2(A,k) \ar[d]^{\tau_{\ul}^*}\\
((t)/(t^2))^\# \ar@{=}[r]& \coh^1(k[t],k)\ar[r]^>>>>>>>>{\sim} & \coh^2(k[t]/(t^{\ell}),k)
}}
\end{xy}
\end{equation}
The leftmost column comes from the isomorphism 
$\coh^1(k[X_1, \dots, X_m],k) \cong \Hom_{{\rm alg}}(k[X_1, \dots, X_m],k) 
= \Hom_k( I/I^2, k) = (I/I^2)^\#$. 
\sloppy
{

}

By construction of the diagram (\ref{lambdadual}), 
the leftmost vertical map is the dual to the map induced by  
$ t \mapsto \lambda_1^{\ell}X_1 + \dots +\lambda_n^{\ell}X_n$. Thus, it sends $z_i$ to $\lambda_i^\ell z$, 
where $z$ is the dual basis to $t$ in $((t)/(t^2))^\#$. 
Therefore,  the rightmost vertical map sends $y_i$ to $\lambda_i^\ell y$.  Since at least one of $\lambda_i$ is nonzero, 
we conclude that $\tau_{\ul}^*$ is surjective onto $\coh^2(k[t]/(t^\ell),k)$, and thus onto $\coh^{ev}(k[t]/(t^\ell),k)$.

\vspace{0.1in}
(ii)
Let $\lambda ' = g\cdot \lambda$.
By the definition (\ref{gaction}) of this action, as $q^{\ell}=1$,
we have $(\lambda_i')^{\ell}=\lambda_i^{\ell}$ for all $i$.
It now follows from the proof of (i) that
$\tau_{\ul}^*=\tau_{\ul'}^*$.
\end{proof}

The lemma implies that we can define a map
\begin{equation}\label{eqn:defnpsi}
\Psi: V^r_A \to V^c_A
\end{equation}
by  sending $\ul\in V_A^r/G ={\mathbb P}^{m-1}/G$ 
to the homogeneous prime ideal $\Ker (\tau_{\ul}^*)$ of $\coh^*(A,k)$. 
The following proposition is an immediate consequence of Lemma \ref{lem:well}.

\begin{prop}
\label{prop:homeo}
$
\Psi: V^r_A \to V^c_A
$
is a homeomorphism.
\end{prop}

\begin{proof}
As $V^r_A\cong {{\mathbb P}}^{m-1}/G$, first define 
$\widetilde \Psi: {\mathbb P}^{m-1} \to V_A^c$
by 
$$\widetilde {\Psi}(\ul)  = \Ker (\tau_{\ul}^*).$$

As it is shown in the proof of Lemma \ref{lem:well}(i),  $\tau_{\ul}^*(y_i) = 
\lambda^{\ell}_i y$, where $y$ is the degree $2$ generator of $\coh^*(k[t]/(t^\ell),k)$.  
Thus, $\Ker (\tau_{\ul}^*)$ is generated 
by the elements $\sum a_iy_i \in \coh^2(A,k)$ such that $\sum a_i \lambda^{\ell}_i = 0$. 
We get 
$$
\widetilde \Psi(\ul) = [\lambda_1^{\ell} : \lambda_2^{\ell}: \cdots 
:\lambda_m^{\ell}].
$$

Finally, since $[\lambda_1^{\ell} : \lambda_2^{\ell}: \dots :\lambda_m^{\ell}] = 
[\mu_1^{\ell} : \mu_2^{\ell}: \dots :\mu_m^{\ell}]$ if and only if there exists $g \in G$ such that
$[\lambda_1: \dots : \lambda_m] = g \cdot [\mu_1: \dots : \mu_m]$, we conclude that $\widetilde \Psi: {\mathbb P}^{m-1} \to V_A^c$ factors through
$\Psi: V_A^r\to V_A^c$ and, moreover, that $\Psi$ is a homeomorphism.  
\end{proof}

\begin{rem}
The rank variety $V_A^r$
can be identified with
$\Proj S^*((I/I^2)^\#)^G \cong  (\Proj I/I^2)/G = \mathbb P^{m-1}/G$.
The map $\Psi$ 
is then given by the algebraic map 
\begin{equation}
\label{phi}
\psi: \coh^*(A,k) \to   (S^*((I/I^2)^\#))^G,
\end{equation}
$$
y_i \mapsto z_i^\ell.
$$
\end{rem}

We will need the following observation, which is well-known in the case of cyclic finite groups (\cite[3.2]{Ev}).
Let $y$ be the degree $2$ generator of $\coh^*(k[t]/(t^{\ell}),k)$. Then multiplication by $y$ induces an isomorphism 
\begin{equation}\label{lem:period}
\cdot y: \coh^{n}(k[t]/(t^{\ell}),N)\to \coh^{n+2}(k[t]/(t^{\ell}),N)
\end{equation}
for any $n>0$ and any $k[t]/(t^{\ell})$-module $N$.
To see this, note that by the periodicity of the trivial module $k$ arising
from (\ref{eqn:tres}), we have $\coh^{n}(k[t]/(t^{\ell}),N)\cong
\coh^{n+2}(k[t]/(t^{\ell}),N)$ for all $n>0$.
The element $y\in \coh^2(k[t]/(t^{\ell}),k)$ corresponds to the identity map in
$\underline{\Hom}_{k[t]/(t^{\ell})}(k,k)\cong
\underline{\Hom}_{k[t]/(t^{\ell})}(\Omega^2(k),k)$, and its cup product with an
element in $\coh^{n}(k[t]/(t^{\ell}),N)\cong\underline{\Hom}
_{k[t]/(t^{\ell})}(\Omega^n(k),N)$ 
corresponds to composition with the identity 
map from $\Omega^n(k)$ to $\Omega^n(k)$.

Now we are able to determine the support varieties of the modules $V(\ul)$,
which will be used to obtain a connection between the rank and support
varieties of an arbitrary finitely generated $A$-module.

\begin{lem}\label{lem:c}
$V_A^c(V(\underline{\lambda}))= \Psi(\ul)$.
\end{lem}

\begin{proof}
Let $z \in I(k,V(\underline{\lambda}))$ be a homogeneous element of even degree in the annihilator of $\Ext^*_A(k,V(\underline{\lambda}))$ in $\Ext_A^*(k,k)$. 
Since the restriction map  
$$\tau_{\ul}^*:  \Ext^*_A(k,V(\underline{
\lambda}))\rightarrow \Ext^*_{k[t]/(t^{\ell})}
(k,V(\underline{\lambda}))$$ 
is injective by Lemma \ref{lem:a}, we can choose 
$v \in \Ext^*_A(k,V(\underline{\lambda}))$ such that $\tau^*_{\ul}(v) \not = 0$.   Since $z \in I(k,V(\underline{\lambda}))$, 
we conclude that 
$$\tau^*_{\ul}(z)\tau^*_{\ul}(v) =\tau^*_{\ul}(zv) = 0.$$  
Due to the isomorphism (\ref{lem:period})
and as $\tau^*_{\ul}(v) \not = 0$, 
 we get $\tau^*_{\ul}(z) = 0$.
Thus, $z \in \Ker(\tau_{\ul}^*)$.   Since all elements of $I(k,V(\ul))$ of odd 
degree are nilpotent and $\Ker(\tau_{\ul}^*)$ is a prime ideal, 
it follows that 
$$
I(k,V(\underline{\lambda})) \subset \Ker(\tau_{\ul}^*).
$$
This implies that $\Psi(\underline{\lambda})\in V^c_A(k,V(
\underline{\lambda}))$. By Proposition \ref{prop:prop}(iv),
$V^c_A(k,V(\underline{\lambda}))\subseteq
V^c_A(V(\underline{\lambda}))$. Therefore $\Psi (\underline{\lambda})
\in V^c_A(V(\underline{\lambda}))$.

It remains to prove that $V^c_A(V(\ul)) \subset \Psi(\ul)$. Applying 
Proposition \ref{prop:prop}(iv) again, 
it suffices to show $V^c_A(S, V(\ul)) \subset \Psi(\ul)$ for any simple A-module $S$. 
This, in turn, will follow from the inclusion $\Ker (\tau_{\ul}^*) \subset I(S, V(\ul))$.

Let $S$ be a simple $A$-module. We claim that the following diagram commutes: 
\begin{equation}
\label{S}
\begin{xy}*!C\xybox{
\xymatrix{
\Ext^*_A(k,k) \ar[rr]^{-\ot S}\ar[drr]^{\tau_{\ul}^*}
&& \Ext^*_A(S,S) \ar[d]^{\tau_{\ul}^*}\\
&&\Ext_{k[t]/(t^{\ell})}^{*}(k,k) 
}}
\end{xy}
\end{equation}
Indeed, suppose $S=S_{\chi}$ where $\chi: G\rightarrow k^{\times}$ is a character,
so that each $g\in G$ acts as multiplication by $\chi(g)$ and each $X_i$ acts as $0$.
Under the map $ - \otimes S$, an $n$-extension $k\rightarrow M_{n-1}\rightarrow
\cdots\rightarrow M_0 \rightarrow k$ is sent to $S\rightarrow M_{n-1}\otimes S\rightarrow
\cdots \rightarrow M_0\otimes S\rightarrow S$.
The action of $\tau_{\ul}(t)$ on each module $M_i\otimes S$ is as 
$$\Delta(\tau_{\ul}(t))=\sum_{i=1}^m\lambda_i(X_ih_i\otimes h_i + g_ih_i\otimes X_ih_i).$$ 
Since $X_ih_i$ acts by $0$ on $S$, this is the same as the action of
$\sum_{i=1}^m \lambda_i X_i h_i\otimes h_i$, which is $\sum_{i=1}^m \chi(h_i)\lambda_i
X_ih_i\otimes 1$.
Thus, when restricted to $k[t]/(t^{\ell})$ via $\tau_{\ul}^*$, the $n$-extension
is equivalent to $k\rightarrow M'_{n-1}\rightarrow\cdots\rightarrow M'_0\rightarrow k$,
where $M'_i=M_i$ as a vector space, and $t$ acts on $M_i'$ as $\sum_{i=1}^m\chi(h_i)
\lambda_iX_ih_i = \tau_{\ul '}(t)$ with $\ul ' = [\chi(h_1)\lambda_1 : \cdots :
\chi(h_m)\lambda_m]$. Therefore there is a $g\in G$ with $\ul ' = g\cdot \ul$.
By Lemma \ref{lem:well}(ii), $\tau_{\ul}^* = \tau_{\ul '}^*$, and so
the diagram commutes.

Note that the map $ - \otimes S$ in fact identifies $\Ext^*_A(k,k)$ and $\Ext^*_A(S,S)$
as graded vector spaces: An inverse map is given by $- \otimes S^{\#}$ since 
$S\otimes S^{\#}\cong k$.

Consider the following commutative diagram where the vertical arrows are restrictions via $\tau_{\ul}$ and horizontal arrows are actions via Yoneda product:
\begin{equation}
\label{diag}
\begin{xy}*!C\xybox{
\xymatrix{
\Ext^*_A(S,S) \ar[r] \ar[d]^{(\tau_{\ul}^*,\tau_{\ul}^*)} 
\times \Ext^*_A(S,V(\underline{\lambda}))
& \Ext^*_A(S,V(\underline{\lambda})) \ar[d]^{\tau_{\ul}^*}\\
\Ext_{k[t]/(t^{\ell})}^{*}(k,k) \times \Ext_{k[t]/(t^{\ell})}^*(k,V(\underline{\lambda}))
\ar[r]
&\Ext_{k[t]/(t^{\ell})}^*(k,V(\underline{\lambda}))
}}
\end{xy}
\end{equation}
The action  of $\Ext^*_A(k,k)$ on $\Ext^*_A(S,V(\underline{
\lambda}))$ factors through the action of $\Ext^*_A(S,S)$ via $-\otimes S: \Ext^*_A(k,k) \to \Ext^*_A(S,S)$. 
By Lemma \ref{lem:a}, the rightmost vertical arrow of (\ref{diag}) is injective.
Let $\alpha\in \Ker (\tau_{\ul}^*)$ in $\Ext^*_A(k,k)$. 
Commutativity of (\ref{S}) implies that 
$\alpha\ot S$ is also in the kernel 
of $\tau_{\ul}^*: \Ext^*_A(S,S) \to \Ext^*_{k[t]/(t^{\ell})}(k,k)$. This means that for every $\beta\in\Ext^*_A(S,V(\underline{\lambda}))$,
we have $\tau_{\ul}^*((\alpha\ot S)\cdot \beta)=0$.
As this $\tau_{\ul}^*$ is injective, this implies $(\alpha\ot S)\cdot\beta=0$,
that is $\alpha\in I(S,V(\underline{\lambda}))$.
Since this holds for any $\alpha \in \Ker (\tau_{\ul}^*)$, we get $\Ker (\tau_{\ul}^*) \subset I(S, V(\ul))$ as required.  
\end{proof}

Finally we use the modules $V(\ul)$ and $L_{\zeta}$ to prove equivalence
of the rank and support varieties. 

\begin{thm}\label{thm:equiv}
Let $M$ be a finitely generated $A$-module. Then
$$\Psi(V^r_A(M))= V^c_A(M).$$
\end{thm}

\begin{proof}
Let $\underline{\lambda}\in V^r_A(M)$, where we abuse notation
by identifying an element $\ul\in {\mathbb P}^{m-1}$ with its $G$-orbit.
Then
$\underline{\Hom}_A(V(\underline{\lambda}),M)\neq 0$ by 
Lemma \ref{lem:projective}(i).
By periodicity (Lemma \ref{lem:vlambda}(v)), 
\begin{equation*}
\Ext^{2n}_A(V(
\underline{\lambda}),M)=\underline{\Hom}_A(V(\underline{\lambda}),M)
\neq 0
\end{equation*}
 for every positive integer $n$. 
This implies that $V_A^c(V(\ul),M) \not = \emptyset$. 
Indeed, suppose $V_A^c(V(\ul), M) = \emptyset$.  
Then $\sqrt{I(V(\ul), M)} = \coh^{>0}(A,k)$.  
Note that $\Ext^*_A(V(\ul),M)\cong \coh^*(A,V(\ul)^\# \otimes M)$ 
(see the appendix), 
and this is a finitely generated module 
over $\coh^*(A,k)$ by Lemma \ref{finite}.
This implies that $\Ext_A^n(V(\ul),M) = 0$ for all
sufficiently large $n$, contradicting what we found above. 
Thus, $V_A^c(V(\ul), M) \not = \emptyset$. 
\sloppy 
{

}

Proposition \ref{prop:prop}(iii) and Lemma \ref{lem:c} imply
that $V_A^c(V(\ul), M) \subset V_A^c(V(\ul)) \cap V_A^c(M) = \Psi(\ul) \cap \va(M)$. 
Since $V_A^c(V(\ul), M) \not = \emptyset$, we conclude that $\Psi(\ul) \subset \va(M)$. 
This proves the containment $\Psi(V^r_A(M))\subseteq V^c_A(M)$.

Now suppose $\Psi(\ul) \in \va(M)$.  
We can find a finite set of homogeneous elements $\zeta_1, \dots, \zeta_{m-1}$
in $\coh^*(A,k)$
(where $m-1$ is the projective dimension of $V^c_A$)
such that  $\Psi(\ul) = \bigcap \langle \zeta_i \rangle$.
Letting $L_{\ul} = L_{\zeta_1} \otimes \cdots \otimes 
L_{\zeta_n}$ (see (\ref{exact})), 
Corollary \ref{cor:realization} implies that $\va(M \otimes L_{\ul}) = \Psi(\ul)$. 
By the first part of the proof, it follows that
$V_A^r(M \otimes L_{\ul}) \subset {\ul}$.
Since $\va(M\otimes L_{\ul}) = \Psi(\ul) \not = \emptyset$, Proposition \ref{prop:prop}(v) implies that
$M \otimes L_{\ul}$ is not  projective.
Theorem \ref{thm:dade} implies that 
$V_A^r(M \otimes L_{\ul}) \not = \emptyset$. Hence, $V_A^r(M \otimes L_{\ul}) = \ul$.

For each $i$, we have a
short exact sequence $0 \to M \otimes L_{\zeta_i} \to M \otimes \Omega^n (k) \to M \to 0$  obtained by applying $M \otimes -$ 
to (\ref{exact}).
By Proposition \ref{prop:rank}(iv), 
$V_A^r(M \otimes L_{\zeta_i}) \subset V_A^r(M) \cup V_A^r(M \otimes \Omega^n (k) )$. 
Since $M \otimes \Omega^n (k)  \cong \Omega^n (M)$ in the stable module category,
Proposition \ref{prop:rank}(iii) implies that 
$V_A^r(M \otimes L_{\zeta_i} ) \subset V_A^r(M)$. Proceeding by induction, we conclude that
$$
V_A^r(M \otimes L_{\ul}) \subset V_A^r(M).$$
Thus, ${\ul} \in V_A^r(M)$.  Since this holds for any $\ul$ such that $\Psi(\ul) \in \va(M)$,  we get the desired inclusion
 $V_A^c(M) \subseteq \Psi(V_A^r(M))$.
\end{proof}

As a consequence of Lemma \ref{lem:c} and Theorem \ref{thm:equiv},
we may now record the rank variety of $V(\ul)$.

\begin{cor}\label{rankvlambda}
$V^r_A(V(\underline{\lambda}))=\{\underline{\lambda}\}$.
\end{cor}


\section{Varieties for modules of truncated polynomial algebras}
\label{sec:trunc}

Our results have consequences for modules of the truncated polynomial
algebra $\Lambda = k[X_1,\ldots,X_m]/(X_1^{\ell},\ldots,X_m^{\ell})$,
which we give in this section.
We define the rank variety of a $\Lambda$-module $M$ by
\begin{equation}\label{eqn:ranklambda}
  V^r_{\Lambda}(M)= V^r_A(M\uparrow ^A),
\end{equation}
where the rank variety of the induced
$A$-module $M\!\uparrow^A = A\ot_{\Lambda}M$ is given in Definition
\ref{defn:rank}.
Since $A$ is free as a $\Lambda$-module, induction from $\Lambda$ to $A$
is well-behaved.   
The rank variety of the trivial $\Lambda$-module is $\mathbb P^{m-1}$
since $\tau(\ul)$ acts trivially on $A\otimes_{\Lambda}k$, for any $\ul$.
We will use the notation $V^r_{\Lambda} 
=V^r_\Lambda(k) = V^r_A(k\uparrow ^A) \cong {\mathbb P}^{m-1}/G$.

\vspace{0.1in}

\begin{rems}\label{crit}
(i) An alternative definition of rank varieties for $\Lambda$-modules
is given in \cite[Rk.\ 4.7(2)]{BEH}. We expect that our definition is
equivalent to this one, however we only have a proof that they are equivalent
in case $\ell = 2$ (see Remark \ref{rem:EH} below).

(ii) Viewing the $A$-modules $V(\ul)=A\cdot \tau_{\ul}(t)^{\ell -1}$ 
(see Section 2) as $\Lambda$-modules by restriction, we have the 
following characterization of the rank variety for a $\Lambda$-module $M$:
$V^r_{\Lambda}(M)$ consists of all $\ul\in {\mathbb{P}}^{m-1}/G$
such that $\underline{\Hom}_{\Lambda}(V(\ul),M)\neq 0$
(cf.\ \cite[Lem.\ 3.7(2) and Defn.\ 4.1]{EH1}).
This is a consequence of Lemma \ref{lem:projective}(i) since the
Eckmann-Shapiro Lemma 
implies that $\underline{\Hom}_A(V(\ul),M\!\uparrow^A)
\cong \underline{\Hom}_{\Lambda}(V(\ul),M)$ as a result of an isomorphism
between induced and coinduced modules. The isomorphism
$\Hom_{\Lambda}(A,M)\cong M\!\uparrow ^A$ 
is given by sending $f\in\Hom_{\Lambda}(A,M)$ to $\sum_{g\in G}g\ot f(g^{-1})$.
The left $A$-module structure of $\Hom_{\Lambda}
(A,M)$ is standard, given by $(a\cdot f)(b) = f(ba)$
for all $a,b\in A$, $f\in\Hom_{\Lambda}(A,M)$. 
\end{rems}

\smallskip

We define the support variety $V^c_\Lambda(M)$ using the action of 
$\coh^*(\Lambda, k)=\Ext^*_{\Lambda}(k,k)$ on $\coh^*(\Lambda, M)=
\Ext^*_{\Lambda}(k,M)$ by Yoneda composition.
Note that $\coh^*(\Lambda,k)$ is the tensor product of a symmetric
algebra with an exterior algebra, each in $m$ variables; this may be
seen by using (\ref{eqn:tres}) and the K\"unneth formula.
In particular, $\coh^*(\Lambda,k)_{red} \cong \coh^*(A,k)$.
Let $\coh^{\bu}(\Lambda,k)$ denote the subalgebra of $\coh^*(\Lambda,k)$
that is the sum of all even degree components when ${\rm char} \ k \neq 2$,
and $\coh^{\bu}(\Lambda,k)=\coh^*(\Lambda,k)$ when ${\rm char} \ k = 2$.

\begin{defn} 
\label{lambda-support} Let $V^c_\Lambda = \Spec \coh^\bu(\Lambda,k)$. 
For a finitely generated $\Lambda$-module $M$, 
define $V^c_\Lambda(M)=V^c_\Lambda(k,M)$ to be the closed subset of 
$V^c_\Lambda$ defined by the annihilator ideal 
$\Ann_{\coh^\bu(\Lambda,k)}\coh^*(\Lambda,M)$.
 \end{defn}

\begin{rem}
Since $\Lambda$ does not have a Hopf algebra structure, we cannot consider the action of $\coh^\bu(\Lambda,k)$ on $\Ext_{\Lambda}^*(M,M)$ 
as it is usually done to define support varieties when the Hopf algebra structure is present.  Nonetheless, since $\Lambda$ 
has a unique trivial module $k$, Definition \ref{lambda-support} is parallel to the standard definition in view of Proposition \ref{prop:prop}(iv).
\end{rem}

The map $\Psi$ of (\ref{eqn:defnpsi}) can be identified with a map
$\Psi : V^r_\Lambda \to  V^c_\Lambda$ via 
$V^r_\Lambda = V^r_A$ and $V^c_\Lambda  = V_A^c$; the second identification 
comes from the isomorphism $\coh^*(\Lambda, k)_{red} \cong \coh^*(A,k)$. 
We will show that $\Psi$ takes the rank variety of a finitely generated
$\Lambda$-module to its support variety. 

\begin{thm}\label{J:equiv}
Let $M$ be a finitely generated $\Lambda$-module.
Then 
$$ \Psi(V^r_{\Lambda}(M)) = V^c_{\Lambda}(M).$$
\end{thm}

\begin{proof}
By the definition (\ref{eqn:ranklambda}) of $V^r_{\Lambda}(M)$ and Theorem
\ref{thm:equiv}, we need only
check that $V^c_\Lambda(M) = V_A^c(M\uparrow^A)$. 
By Proposition \ref{prop:prop}(iv) it suffices to show that 
$V^c_\Lambda(M) = V_A^c(S,M\uparrow^A)$ for any simple $A$-module $S$. 
By arguments similar to those in Remark \ref{crit}(ii),
$\Ext_A^*(S, M\uparrow^A) \cong \Ext_\Lambda^*(S, M) \cong \Ext_\Lambda^*(k, M)$.   This  is an isomorphism of $\coh^*(A,k)$-modules where $\coh^*(A,k)$ acts on 
$\Ext^*_{\Lambda}(k,M)$ via the embedding 
$\coh^*(A,k) \hookrightarrow \coh^*(\Lambda,k)$.  
Since $\coh^*(A,k) = \coh^\bu(\Lambda,k)_{red}$, 
the variety of the annihilator of  $\Ext_\Lambda^*(k, M)$ is determined by the action of $\coh^*(A,k)$. The statement follows.
\end{proof}

We now explain the connection between our results and the work on support varieties which was done from the point of view of 
Hochschild cohomology in \cite{EH2}, \cite{EHSST}, \cite{SS}.

We will show that the rank variety of a $\Lambda$-module is also equivalent
to its {\em Hochschild} support variety defined 
as follows via a particular choice of subalgebra of
the Hochschild cohomology ring
$\HH^*(\Lambda)=\Ext^*_{\Lambda^e}(\Lambda,\Lambda)$, where
$\Lambda^e=\Lambda\ot\Lambda^{op}$.
We use the definition of (Hochschild) support varieties for modules of finite dimensional
algebras given in \cite{SS} and developed further in \cite{EHSST}.
Let 
$$H=\coh^*(A,k).$$
We will show that $H$ embeds canonically as a subalgebra of $\HH^*(\Lambda)$.
To do this, we will need to consider the following subalgebra of 
$A^e=A\ot A^{op}$,
$${\mathcal D} =\Lambda^e\rtimes\delta(G) =\bigoplus_{g\in G}(\Lambda g
  \otimes \Lambda g^{-1})\subset A^e,$$
where $\delta(G)=\{(g,g^{-1})\mid g\in G\}\cong G$ acts on $\Lambda^e$ by
the left and right actions induced by the action of $G$ on $\Lambda$.
Note that $\mathcal D$ contains the subalgebra $\delta(A)\cong A$,
where $\delta(a)=\sum a_1\ot S(a_2)$ (see 
Lemma \ref{induction} below). As noted in the proof of
Lemma \ref{ae}, $A^e$ is projective as a right $\delta(A)$-module
under multiplication. Note that $A^e=\oplus_{g\in G} ((g\ot 1){\mathcal D})$
as a right $\delta(A)$-module, that is $\mathcal D$ is a direct summand
of the projective $\delta(A)$-module $A^e$, and so is projective itself.

Let $P_{\bu}\rightarrow k$ be an $A$-projective resolution of the trivial
$A$-module $k$. The isomorphism $A\cong A^e\ot_{\delta(A)}k$
of $A^e$-modules given in Lemma \ref{induction} restricts to an isomorphism
$\Lambda\cong {\mathcal D}\ot_{\delta(A)} k$ of $\mathcal D$-modules.
Therefore induction to $\mathcal D$ yields a $\mathcal D$-projective
resolution of the $\mathcal D$-module $\Lambda$:
$$
  \cdots \rightarrow {\mathcal D}\ot_{\delta(A)} P_1 \rightarrow
  {\mathcal D} \ot_{\delta(A)} P_0\rightarrow \Lambda\rightarrow 0.
$$
Induction further to $A^e$ yields an $A^e$-projective resolution
$$
  \cdots \rightarrow A^e\ot_{\delta(A)} P_1 \rightarrow
  A^e \ot_{\delta(A)} P_0\rightarrow A\rightarrow 0.
$$
Now suppose $\xi\in\coh^n(A,k)$, and identify $\xi$ with a representative $A$-map
$\xi :P_n\rightarrow k$.
Induction yields a $\mathcal D$-map 
$\xi ' : {\mathcal D}\ot_{\delta(A)} P_n\rightarrow
\Lambda$, and further induction
yields an $A^e$-map 
$\xi '':A^e\ot_{\delta(A)} P_n\rightarrow A$.
The induction from $\delta(A)$ to $A^e$ results in precisely the embedding of
$\coh^n(A,k)$ into $\HH^n(A)$ given in Lemma \ref{ae}.
Therefore, the map sending $\xi$ to $\xi '$ is an embedding of
$\coh^n(A,k)$ into $\Ext^n_{\mathcal D}(\Lambda,\Lambda)$.
Note that $\Ext^n_{\mathcal D}(\Lambda,\Lambda)
\cong \Ext^n_{\Lambda^e}(\Lambda,\Lambda)^G =\HH^n(\Lambda)^G$
since the characteristic of $k$ is relatively prime to the order of 
$G$.
This provides the embedding of $H=\coh^*(A,k)\hookrightarrow
\HH^*(\Lambda)^G\hookrightarrow \HH^*(\Lambda)$.

\begin{rem} The embedding $H=\coh^*(A,k)\hookrightarrow \HH^*(\Lambda)$ can be described explicitly as follows.
In case $m=1$, identify $\Lambda$ with $k[t]/(t^{\ell})$. There is a periodic
$\Lambda^e$-free resolution of $\Lambda$:
$$
  \cdots\stackrel{v\cdot}{\longrightarrow} \Lambda^e 
        \stackrel{u\cdot}{\longrightarrow} \Lambda^e
\stackrel{v\cdot}{\longrightarrow} \Lambda^e 
        \stackrel{u\cdot}{\longrightarrow} \Lambda^e
  \stackrel{m}{\longrightarrow} \Lambda \rightarrow 0
$$
where $u=t\ot 1 - 1\ot t$ and $v=t^{\ell -1}\ot 1 + t^{\ell -2}\ot t
+\cdots + 1\ot t^{\ell -1}$ \cite[Exer.\ 9.1.4]{weibel94}.
Using this resolution, one computes $\HH^n(\Lambda)\cong\Lambda/(t^{\ell -1})$ and
$\HH^n(\Lambda)^G \cong k$ for all $n> 0$.
As $\coh^n(A,k)=k$ for all even $n\geq 0$, and $0$ for all
odd $n>0$, the embedding $\coh^*(A,k)\hookrightarrow \HH^*(\Lambda)^G$ is forced
to be an isomorphism onto $\HH^*(\Lambda)^G_{red}$ in case the
characteristic is not $2$.
Apply the K\"{u}nneth Theorem to obtain the embedding for all $m$.
\end{rem}

We must verify that $\Lambda$ and $H$ satisfy the properties required by the theory
of support varieties defined via Hochschild cohomology in \cite{EHSST}:
As $\Lambda$ is local, it is an indecomposable algebra. 
The cohomology algebra $H=\coh^*(A,k)$ is a polynomial ring in $m$ 
variables, each of degree 2  (\ref{eqn:cohA}), 
so it is a commutative, noetherian, graded
subalgebra of $\HH^*(\Lambda)$. The assumption $H^0=Z(\Lambda)=\Lambda$
does not hold. However the generators of $\Lambda$ are nilpotent and so
we may consider $H$ to be the reduced version
of the subalgebra of $\HH^*(\Lambda)$ generated by $H$ and $\Lambda$,
and for the purpose of defining varieties it suffices just to consider $H$.
Thus $H$ essentially satisfies assumption (Fg1) of \cite{EHSST}.
Since $\coh^*(A,k) \cong\Ext^*_{\Lambda}(k,k)^G$,  we get that  $\Ext^*_{\Lambda}(k,k)$ is a 
finitely generated module over $H$,
however we must check that the usual action agrees with that defined via
the subalgebra $\mathcal D$ above. To see this, it suffices to check
that $({\mathcal D}\ot_{\delta(A)}M)\ot_{\Lambda}k\cong M$ as
$\Lambda$-modules, for any $A$-module $M$.
Using the techniques of Lemma \ref{induction}, we may write an arbitrary
element of $({\mathcal D}\ot_{\delta(A)}M)\ot_{\Lambda}k$ as a linear
combination of certain elements
\begin{eqnarray*}
  a\ot b\ot m\ot 1 &=& \sum a_1\ot\varepsilon(a_2)b\ot m\ot 1\\
  &=& \sum a_1\ot S(a_2)a_3b\ot m \ot 1 \\
  &=& \sum 1\ot a_2b\ot a_1m\ot 1\\
  &=& \sum 1\ot 1\ot a_1m\ot \varepsilon(a_2b) = 1\ot 1\ot\varepsilon(b)
     am\ot 1
\end{eqnarray*}
($a,b\in A$, $m\in M$)
since $\Lambda$ is a left coideal subalgebra of $\Lambda$ (that is $a_2b\in
\Lambda$ when $a\ot b\in {\mathcal D}$).
The desired isomorphism $({\mathcal D}\ot_{\delta(A)} M)\ot_{\Lambda}k
\rightarrow M$ is thus given by $a\ot b\ot m\ot 1\mapsto \varepsilon(b)am$,
with inverse $m\mapsto 1\ot 1\ot m\ot 1$. This proves that
(Fg2) of \cite{EHSST} is satisfied.

As in \cite{EHSST,SS} we define the (Hochschild) support variety of a $\Lambda$-module
$M$, with respect to $H=\coh^*(A,k)$, as
\begin{equation}\label{eqn:supplambda}
  V^H_{\Lambda}(M)= \Proj H/\Ann_H \Ext^*_{\Lambda} (M,M)
\end{equation}
where the action of $H$ on $\Ext^*_{\Lambda}(M,M)$ is by $-\ot_{\Lambda}M$
(under the identification of elements of $H$ with $\Lambda^e$-extensions
of $\Lambda$ by $\Lambda$) followed by Yoneda composition. 
We next show that $V^H_{\Lambda}(M)$ is homeomorphic to the support
variety $V^c_{\Lambda}(M)$ given in Definition \ref{lambda-support}.

\begin{thm}\label{thm:equivsupp}
Let $M$ be a finitely generated $\Lambda$-module. Then 
$$V^H_{\Lambda}(M)\cong V^c_\Lambda(M).$$
\end{thm}

\begin{proof}
It was shown in the proof of Theorem \ref{J:equiv} that $V^c_{\Lambda}(M)
= V^c_A(M\uparrow ^A)$.
We will show here that $V^c_A(M\uparrow^A)\cong V^H_{\Lambda}(M)$.
We analyze the following diagram:
\begin{equation*}
\begin{xy}*!C\xybox{
\xymatrix{
  \coh^{*}(A,k) \ar[drr]^{-\ot_k M\uparrow^A}\ar[d]_{A^e\ot_A - }
  &&\\ 
  \HH^{*}(A) \ar[rr]^{-\ot_{A} M\uparrow^A}
  &&\Ext^{*}_A(M\!\uparrow^A,M\!\uparrow^A)\\
  \Ext^{*}_{\mathcal D}(\Lambda,A)\ar[u]^{\cong} 
 && \Ext^{*}_{\Lambda}(M,M\!\uparrow^A)\ar[u]^{\cong}\\
  \HH^{*}(\Lambda)^G \ar[rr]^ 
   {-\ot_{\Lambda} M}\ar[u]
   &&\Ext^{*}_{\Lambda}(M,M)\ar[u]
}}\end{xy}
\end{equation*}
The lower left arrow going up is defined by applying the isomorphism
$\Ext^*_{\mathcal D}(\Lambda,\Lambda)\cong \HH^*(
\Lambda)^G$ discussed above and identifying $\Lambda$ with the
$\mathcal D$-submodule $\Lambda\rtimes 1$ of $A=\Lambda\rtimes G$.
The lower right arrow going up is defined by identifying $M$ with the 
$\Lambda$-submodule $1\otimes M$ of $M\!\uparrow^A$.
The embedding of $\coh^*(A,k)=\Ext^*_A(k,k)$ into $\HH^*(A)=\Ext^*_{A^e}(A,A)$ in the diagram
identifies it with a subalgebra of the image of $\HH^*(
\Lambda)^G$ in $\HH^*(A)$ as a result of the discussions above
on the algebra $\mathcal D$.

The top triangle commutes by Lemma \ref{lem:H}.
This implies that the annihilators of 
$\Ext^{*}_A(M\uparrow^A,M\uparrow^A)$ in $\Ext^{*}_A(k,k)$ and
in the subalgebra $H$ of $\Ext^{*}_{A^e}(A,A)$ coincide.
We have also seen that $H$ may be identified with a subalgebra of
$\HH^{*}(\Lambda)^G$, compatible with these
embeddings. Therefore it remains to check commutativity of the bottom
part of the diagram, since the vertical arrows are injections, implying
the appropriate annihilators will coincide. 
This may be checked directly at the chain level, using the 
identification $\HH^*(\Lambda)^G\cong \Ext^*_{\mathcal D}(\Lambda,\Lambda)$
given above.
\end{proof}

We have an immediate consequence of Theorems \ref{J:equiv} and \ref{thm:equivsupp}.

\begin{cor}\label{cor:equivlambda}
Let $M$ be a finitely generated $\Lambda$-module.
Then
$$
  \Psi(V^r_{\Lambda}(M)) \cong V^H_{\Lambda}(M).
$$
\end{cor}

\begin{rem}\label{rem:EH}
In case $\ell =2$, our support varieties for 
$\Lambda$-modules are equivalent to those defined by Erdmann and Holloway
\cite{EH2} since our choice of $H$ is precisely $\HH^*(\Lambda)$ modulo
nilpotent elements. 
As a consequence of Corollary \ref{cor:equivlambda} and the
analogue of the Avrunin-Scott Theorem in \cite{EH2}, our rank varieties
for $\Lambda$-modules also coincide with those defined in \cite{EH2}.
In case $\ell >2$, we expect our rank varieties for $\Lambda$-modules
to coincide with those defined by Benson, Erdmann, and Holloway
\cite[Rk.\ 4.7(2)]{BEH}, but the cohomological techniques are not
yet available in this case and so perhaps a different approach is needed.
\end{rem}

\section{Appendix}
\label{sec:hochcoh}

\subsection*{Hochschild cohomology of Hopf algebras}
Here we allow $A$ to be any Hopf algebra over the field $k$,
and record some general results.
In particular,
we give connections between the cohomology
$\coh^*(A,k)$ and the Hochschild cohomology $\HH^*(A)=\Ext_{A^e}(A,A)$,
where $A^e=A\ot A^{op}$ acts on $A$ by left and right multiplication.
An embedding of $\coh^*(A,k)$ into $\HH^*(A)$ is deduced by Ginzburg 
and Kumar \cite[Prop.\ 5.6 and Cor.\ 5.6]{GK}. 
Here we give an explicit map expressing such an embedding at the chain
level and the resulting connections between actions on $\Ext^*_A(M,M)$
where $M$ is an $A$-module
(see Lemma \ref{lem:H}).

\begin{lem}\label{induction}
There is an isomorphism of $A^e$-modules 
$A\cong (A^e)\otimes_A k,$
where $A$ is embedded as a subalgebra of $A^e$ via the map
$\delta:A\rightarrow A^e$ defined by $\delta(a)=\sum a_1\ot S(a_2)$.
\end{lem}

\begin{proof}
First note that $\delta$ is indeed injective as $\pi\circ\delta ={\rm id}$
where the linear map $\pi :A^e\rightarrow A$ is defined by
$\pi(a\ot b)=a \ \varepsilon(b)$. As $S$ is an algebra anti-homomorphism,
$\delta$ is an algebra homomorphism. Next define $f:A^e\ot_A k\rightarrow A$
by $f(a\ot b\ot 1)=ab$ and $g:A\rightarrow A^e\ot_A k$ by $g(a)=a\ot 1\ot 1$.
We  check that $f$ and $g$ are inverse maps:
$f\circ g(a)=f(a\ot 1\ot 1)=a$ for all $a\in A$, and
\begin{eqnarray*}
   g\circ f (a\otimes b\otimes 1) &=& ab\otimes 1\otimes 1\\
       &=&\sum ab_1\otimes\varepsilon(b_2)\otimes 1\\
    &=&\sum ab_1\otimes S(b_2)b_3\otimes 1\\
   &=& \sum a\otimes b_2\otimes\varepsilon(b_1) = a\otimes b\otimes 1,
\end{eqnarray*}
for all $a,b\in A$. Similar calculations show that $f$ and $g$ are
both $A^e$-module homomorphisms.
\end{proof}

We will next recall some homological properties of modules for a Hopf
algebra.
For proofs, see \cite[I.\ \S3.1]{Ben}.
Let $U$ be a projective left (respectively, right)
 $A$-module, and $V$ any left (respectively, right) $A$-module.
Then $V\otimes U$ is a projective left (respectively, right) $A$-module.
If $V,W$ are left $A$-modules, then $\Hom_k(V,W)$ is a left $A$-module
under the action
$(a\cdot f)(v) = \sum a_2\cdot f(S(a_1)\cdot v).$
In this way, $\Hom_k(V,W)\cong V^\#\ot W$ as $A$-modules, where $V^\#
=\Hom_k(V,k)$ is an $A$-module similarly, with $k$ taking the trivial action
of $A$.
If $U$, $V$, and $W$ are left $A$-modules,
then there is a natural isomorphism
$\Hom_k(V\otimes U, W)\cong \Hom_k(U,V^\#\ot W)$
of $A$-modules. 
This restricts to a natural isomorphism of vector spaces
$\Hom_A(V\otimes U,W)\cong \Hom_A(U,V^\#\ot W)$, further inducing
an isomorphism of graded vector spaces
$$\Ext^*_A(V\ot U, W)\cong \Ext^*_A(U,V^\#\ot W).$$
Similar results hold if $U,V,W$ are right $A$-modules, where
$\Hom_k(V,W)$ is a right $A$-module under the action 
$(f\cdot a)(v)=\sum f(v\cdot S(a_1))\cdot a_2$.
It is helpful to view $\Hom_A(U,V)$ as 
the subspace of $A$-invariant elements of $\Hom_k(U,V)$,
that is as
$$(\Hom_k(U,V))^A=\{f\in\Hom_k(U,V)\mid a\cdot f=\varepsilon(a)f
\mbox{ for all }a\in A\},$$
equivalent by a straightforward computation \cite[Lem.\ 1]{zhu94}.

We will consider $A$ to be a left $A$-module under the left adjoint action,
that is if $a,b\in A$, 
$$a\cdot b = \sum a_1 bS(a_2).$$
Denote this $A$-module by $A^{ad}$.

\begin{lem}\label{ae}
Assume the antipode $S$ is bijective.
There is an isomorphism of graded algebras $\HH^*(A)\cong \coh^*(A,A^{ad})$.
This induces an embedding of $\coh^*(A,k)$ into $\HH^*(A)$.
\end{lem}

\begin{proof}
As $S$ is bijective, it may be checked
that $S: A\rightarrow A^{op}$ is an isomorphism of right 
$A$-modules, where $A$ acts on the right by multiplication on $A$
and by multiplication by $S(A)$ on $A^{op}$.
This yields an isomorphism of right $A$-modules $A\otimes A\rightarrow
A\otimes A^{op}=A^e$.
Now $A\otimes A$ is projective as it is a tensor product of projective
modules.
Thus $A^e$ is a projective $A$-module, the action of $A$ being
 multiplication by $\delta(A)$.
We may therefore apply
the Eckmann-Shapiro Lemma, together with Lemma \ref{induction}, to obtain
an isomorphism of vector spaces
$\Ext^*_A(k,A^{ad})\cong \Ext^*_{A^e}(A,A)$.
This is in fact an isomorphism of algebras. 
(The correspondence of cup products follows from a generalization of the 
proof of \cite[Prop.\ 3.1]{SW} from group algebras to Hopf algebras. 
See also \cite[\S5.6]{GK}.)
Note that the trivial module $k$ is a direct summand of $A^{ad}$,
with complement the augmentation ideal $\Ker( \varepsilon)$.
This results in identification of $\coh^*(A,k)=\Ext^*_A(k,k)$ with a
subalgebra of $\HH^*(A)=\Ext^*_{A^e}(A,A)$.
By construction, this identification is given explicitly 
by the map $A^e\ot_A -$ on extensions.
\end{proof}

The next lemma gives the connection between actions on $\Ext^*_A(M,M)$
that is used in Section \ref{sec:trunc}.

\begin{lem} \label{lem:H}
Assume the antipode $S$ is bijective.
Let $M$ be a finitely generated left $A$-module. Then the following
diagram commutes:
\begin{equation*}
\begin{xy}*!C\xybox{
\xymatrix{
   \coh^{*}(A,k) 
     \ar[drr]^{-\otimes_k M}\ar[d]_{A^e\ot_A -}\\
  \HH^*(A) \ar[rr]^{-\otimes_A M}
   &&\Ext^*_A(M,M)
}}\end{xy}
\end{equation*}
\end{lem}

\begin{proof}
The leftmost map $A^e\ot_A -$ is induced by the diagonal embedding 
$\delta:A\rightarrow A^e$ defined in Lemma \ref{induction}.
The bottommost map
involves right multiplication of $A^e$ by $A\cong 1\otimes A$.
Note that $-\otimes_A M$ does indeed take $n$-extensions to $n$-extensions:
We may assume all modules in an $n$-extension of $A$ by $A$ as 
$A^e$-modules are free over $A$. Then $-\ot_A M$ takes such an
exact sequence to another exact sequence.

Now suppose $0\rightarrow k\rightarrow P_n\rightarrow \cdots\rightarrow
 P_1\rightarrow k\rightarrow 0$ is an $n$-extension of $A$-modules
representing an element of $\coh^n(A,k)$. For each $i$, define $k$-linear
maps
\begin{eqnarray*}
  f_i &:& (A^e\ot_A P_i)\ot_A M \rightarrow P_i\ot_k M\\
  g_i &:& P_i\ot_kM \rightarrow (A^e\ot_A P_i)\ot_A M
\end{eqnarray*}
by $f_i(a\ot b\ot p \ot m)=\sum a_1p\ot a_2bm$ and
$g_i(p\ot m)=1\ot 1\ot p\ot m\in (A\ot A^{op})\ot_A P_i\ot_A M$.
Then clearly
$f_i\circ g_i=\id$. On the other hand,
\begin{eqnarray*}
  g_i\circ f_i(a\ot b\ot p\ot m) & = & \sum 1\ot 1\ot a_1p\ot a_2bm\\
   &=& \sum a_1\ot S(a_2)\ot p\ot a_3bm\\
   &=&\sum a_1\ot S(a_2)a_3b\ot p\ot m\\
   &=& \sum a_1\ot \varepsilon(a_2)b\ot p\ot m = a\ot b\ot p\ot m.
\end{eqnarray*}
Therefore $f_i$ and $g_i$ are inverse maps.
We check that $g_i$ is an $A$-map ($f_i$ is easier and is left to the 
reader):
\begin{eqnarray*}
  g_i(a\cdot (p\ot m)) &=& g_i(\sum a_1p\ot a_2m)\\
  &=& \sum 1\ot 1\ot a_1p\ot a_2m\\
  &=& \sum a_1\ot S(a_2)\ot p\ot a_3m\\
  &=& \sum a_1\ot S(a_2)a_3\ot p\ot m\\
  &=& a\ot 1\ot p\ot m = a\cdot g_i(p\ot m).
\end{eqnarray*}
It is straightforward to check that the $f_i,g_i$ are chain maps.
Therefore an $n$-extension of $A$-modules $k$ by $k$ is taken to the
same extension of $M$ by $M$, either way around the diagram.
\end{proof}

\subsection*{Computational lemma}
The following result is needed in the proof of Lemma \ref{lem:a}.
Let $A=\Lambda\rtimes G$ as in Section 2, and recall the notation
$\tau_{\ul}(t)=\sum_{i=1}^m \lambda_iX_ih_i$ where
$\ul = (\lambda_1,\ldots,\lambda_m)\in k^m$.
If $S$ is a simple $A$-module, let $e_S$ be the primitive central
idempotent of $kG$ corresponding to $S$.

\begin{lem}
\label{compute} 
Let $a\in A$, $\ul\in k^m\setminus \{0\}$, and $S$ a simple $A$-module. 
If $\tau_{\ul}(t)a=0$ and $a\tau_{\ul}(t)^{\ell -1}$ is a  
scalar multiple of $e_S X_1^{\ell -1}\cdots X_m^{\ell -1}$,
then $a\tau_{\ul}(t)^{\ell -1} = 0$.
\end{lem}

\begin{proof}
Write $a=\sum_{0\leq a_i,b_i\leq \ell -1} \alpha_{a_1,\ldots,a_m,
b_1,\ldots,b_m} X_1^{a_1}\cdots X_m^{a_m} g_1^{b_1}\cdots g_m^{b_m}$ 
for some scalars $\alpha_{a_1,\ldots,a_m,b_1,\ldots,b_m}
 \in k$.
Assume without loss of generality that $\lambda_1=1$.
We use a $q$-multinomial formula for $\tau_{\ul}(t)^{\ell -1}$, which may
be obtained from the $q$-binomial formula (stated in the proof of 
Lemma \ref{lem:tau}) and induction on $m$. We need the notation
$$
  \binom{n}{s_1,\ldots,s_m}_q = \frac{(n)_q!}{(s_1)_q!\cdots (s_m)_q!}.
$$
Assuming $\tau_{\ul}(t) a =0$ and  $\lambda_1=1$,
\begin{eqnarray*}
0\!\!\!&=&\!\!\! \tau_{\ul}(t)^{\ell -1}a\\
  \!\!\!&=&\!\!\! \left(\sum_{\substack {0\leq s_i\leq \ell -1\\s_1+\cdots +s_m =\ell -1}}
  \binom{\ell -1}{s_1,\ldots,s_m}_q \lambda_2^{s_2}\cdots\lambda_m^{s_m}
   X_1^{s_1}\cdots X_m^{s_m}h_1^{s_1}\cdots h_m^{s_m} \right)\cdot\\
&&\hspace{1cm}\left(\sum_{0\leq a_i,b_i\leq \ell -1}\alpha_{a_1,\ldots,a_m,b_1,\ldots,b_m}
   X_1^{a_1}\cdots X_m^{a_m}g_1^{b_1}\cdots g_m^{b_m}\right)\\
\!\!\!&=&\!\!\! \sum_{\substack {0\leq s_i\leq \ell -1\\s_1+\cdots +s_m=\ell -1}}
  \sum_{0\leq a_i,b_i\leq \ell -1} \alpha_{a_1,\ldots,a_m,b_1,\ldots,b_m}
  \binom{\ell -1}{s_1,\ldots,s_m}_q \lambda_2^{s_2}\cdots \lambda_m^{s_m}\cdot\\
  &&\hspace{1cm}q^{a_1s_2+(a_1+a_2)s_3+\cdots +(a_1+\cdots +a_{m-1})s_m}
  X_1^{a_1+s_1}\cdots X_m^{a_m+s_m}g_1^{b_1+s_2+\cdots +s_m}\cdots
  g_m^{b_m} \\
\!\!\!&=&\!\!\!\sum_{0\leq a_i,b_i\leq \ell -1}\!\left(\sum_{\substack {0\leq s_i\leq \ell -1\\
   s_1+\cdots +s_m =\ell -1 }} \!\!\!\alpha_{a_1-s_1,\ldots,a_m-s_m,b_1-s_2-\cdots
  -s_m,\ldots,b_m} \binom{\ell -1}{s_1,\ldots,s_m}_q \lambda_2^{s_2}\cdots
  \lambda_m^{s_m}\cdot \right.\\
  && \left. q^{(a_1-s_1)s_2+(a_1+a_2-s_1-s_2)s_3+\cdots + (a_1+\cdots +a_{m-1}
  -s_1-\cdots -s_{m-1})s_m} \right) X_1^{a_1}\cdots X_m^{a_m} g_1^{b_1}\cdots
  g_m^{b_m},
\end{eqnarray*}
where we define $\alpha_{a_1,\ldots,a_m,b_1,\ldots,b_m}=0$ if any one
of $a_1,\ldots,a_m$ is negative.
Each coefficient is thus $0$, that is
\begin{equation}\label{coeff1}
 0=\sum_{\substack {0\leq s_i\leq \ell -1\\s_1+\cdots +s_m =\ell -1}}
  \alpha_{a_1-s_1,\cdots,a_m-s_m,b_1-s_2-\cdots -s_m, \cdots, b_m}
  \binom{\ell -1}{s_1,\cdots,s_m}_q \lambda_2^{s_2}\cdots \lambda_m^{s_m}\cdot
\end{equation}

\vspace{-.1in}

$$ \hspace{1in} q^{(a_1-s_1)s_2 +(a_1+a_2-s_1-s_2)s_3+\cdots +(a_1+\cdots +a_{m-1} -s_1-\cdots
  -s_{m-1})s_m}.$$
By similar calculations we have
\begin{eqnarray*}
a\tau_{\ul}(t)^{\ell -1}
  \!\!\!&=&\!\!\! \sum_{0\leq a_i,b_i\leq \ell -1}\!\!\left(\sum_
  {\substack {0\leq s_i\leq \ell -1\\
  s_1+\cdots +s_m= \ell -1}} \!\!\!\alpha_{a_1-s_1,\ldots,a_m-s_m, b_1-s_2-\cdots
  -s_m,\ldots,b_m} \binom{\ell -1}{s_1,\ldots,s_m}_q\cdot\right. \\
&&\hspace{.6in} \left.\lambda_2^{s_2}\!\cdots\!
  \lambda_m^{s_m}
  q^{(b_1-s_2-\cdots -s_m)s_1+\cdots + b_ms_m}\right)
   X_1^{a_1}\cdots X_m^{a_m} g_1^{b_1}\cdots g_m^{b_m}.
\end{eqnarray*}
If this were  a nonzero element of $ke_S X_1^{\ell -1}\cdots X_m^{\ell -1}$, 
then for any $m$-tuple $(b_1,\ldots,b_m)$, the coefficient of $X_1^{\ell -1}
\cdots X_m^{\ell -1} g_1^{b_1}\cdots g_m^{b_m}$ would be nonzero, that is
\begin{equation}\label{coeff2}
  0\neq \sum_{\substack {0\leq s_i\leq \ell -1\\s_1+\cdots +s_m=\ell -1}}
 \alpha_{\ell -1-s_1 ,\ldots, \ell -1-s_m, b_1-s_2-\cdots -s_m,\cdots, b_m}
  \binom{\ell -1}{s_1,\ldots,s_m}_q \lambda_2^{s_2} \cdots \lambda_m^{s_m}\cdot
\end{equation}

\vspace{-.1in}

$$\hspace{1in}q^{(b_1-s_2-\cdots -s_m)s_1 +\cdots + b_ms_m}.
$$
Letting $a_1=\cdots =a_m =\ell -1$ in (\ref{coeff1}) and comparing with
(\ref{coeff2}), we claim that we may choose $(b_1,\ldots,b_m)$
so that these coefficients are equal, a contradiction:
We will find a solution to 
$$
  q^{(\ell -1-s_1)s_2 +\cdots + ((m-1)(\ell -1) -s_1-\cdots -s_{m-1})s_m}
  = q^{(b_1-s_2-\cdots -s_m)s_1 +\cdots +b_ms_m}.
$$
This equation is equivalent to
$
  q^{-s_2-2s_3-\cdots -(m-1)s_m} = q^{b_1s_1+\cdots +b_ms_m}
$, which has solution
$b_1 =0, \ b_2=-1, \ldots,  \ b_m = -(m-1)$.
Thus if $a\tau_{\ul}(t)^{\ell -1}$ is a
scalar multiple of $e_S X_1^{\ell -1}\cdots X_m^{\ell -1}$, 
it must be $0$.
\end{proof}

\quad

\end{document}